\documentclass[12 pt]{amsart}

\usepackage{graphicx}
\usepackage{color}
\usepackage{amssymb, amsmath, amsthm}
\usepackage{mathptmx}

\usepackage{epsfig}
\usepackage{epstopdf}
\usepackage{graphicx}

\theoremstyle{plain}

\newtheorem{theorem}{Theorem}
\newtheorem{lemma}{Lemma}

\newtheorem{question}{Question}

\newcommand{\Z}{\mathbb{Z}}

\newcommand{\bi}{\begin{itemize}}
\newcommand{\ei}{\end{itemize}}
\newcommand{\be}{\begin{enumerate}}
\newcommand{\ee}{\end{enumerate}}

\newcommand{\n}{\beta}

\newcommand{\emp}{\emptyset}
\newcommand{\X}{\times}

\newcommand{\eps}{\epsilon}

\newcommand{\A}{\alpha}
\newcommand{\pd}{\partial}

\numberwithin{definition}{section}
\numberwithin{example}{section}
\numberwithin{lemma}{section}
\numberwithin{theorem}{section}
\numberwithin{corollary}{section}

\begin{document}
\title{Bridge spectra of iterated torus knots}
\author{Alexander Zupan}
\thanks{The author is supported by the National Science Foundation under Award No. DMS-1203988.}

\maketitle

\begin{abstract}
We determine the set of all genus $g$ bridge numbers of many iterated torus knots, listing these numbers in a sequence called the bridge spectrum.  In addition, we prove a structural lemma about the decomposition of a strongly irreducible bridge surface induced by cutting along a collection of essential surfaces.
\end{abstract}

\section{Introduction}

Given a knot $K$ in a compact, orientable 3-manifold $M$, a bridge splitting of $(M,K)$ often characterizes topological and geometric aspects of $K$ in $M$.  Defined by Doll \cite{doll} and Morimoto and Sakuma \cite{mosak}, bridge surfaces can be seen as analogues to Heegaard surfaces in 3-manifold theory:  A bridge surface cuts the pair $(M,K)$ into two simple topological pieces, reducing the essence of $(M,K)$ to a gluing map.\\

Every knot has infinitely many bridge surfaces, but we may narrow our search for structure by looking for \emph{irreducible} bridge surfaces, which are not the result of a generic modification to another surface.  Although there is a significant body of research concerning the set of irreducible Heegaard splittings of a 3-manifold, there are few examples of knots $K \subset S^3$ with many irreducible bridge surfaces.  Several classes of knots are known to have unique irreducible bridge spheres, up to unoriented isotopy.  Knots whose non-minimal bridge spheres are reducible are called \emph{destabilizable}, a definition due to Ozawa and Takao \cite{ozawatakao}, and classes of knots known to have this property include
\bi
\item the unknot \cite{otal},
\item 2-bridge knots \cite{otal},
\item torus knots \cite{ozawa},
\item iterated torus knots and iterated cables of 2-bridge knots \cite{zupan}, and
\item more generally, cables of an $mp$-small destabilizable knot \cite{zupan}.
\ei

In contrast, Ozawa and Takao have recently produced the first example of a knot $K \subset S^3$ such that $K$ has two irreducible bridge spheres with different bridge number \cite{ozawatakao}, and Jang has exhibited 3-bridge links with infinitely many distinct 3-bridge spheres \cite{jang}.  In terms of higher genus bridge surfaces, we are aware of only one result:  Scharlemann and Tomova have shown that 2-bridge knots have a unique irreducible bridge surface up to unoriented isotopy \cite{schartom}. \\

Following Doll \cite{doll}, we define the genus $g$ bridge number $b_g(K)$ of a knot $K \subset S^3$:
\[ b_g(K) = \min\{b : K \text{ admits a } (g,b)\text{-splitting}\},\]
\noindent and we introduce the \emph{bridge spectrum} $\mathbf{b}(K)$,
\[ \mathbf{b}(K) = (b_0(K),b_1(K),\dots).\]
A construction called \emph{meridional stabilization} (see Section \ref{pre}) transforms a $(g,b)$-surface into a $(g+1,b-1)$-surface; hence, the bridge spectrum is bounded above by the sequence $(b_0(K),b_0(K)-1,\dots,1,0)$.  In view of this property, we say that $\mathbf{b}(K)$ has a \emph{gap} at index $g$ if $b_g(K) < b_{g-1}(K) - 1$, and it is a simple verification that if $\mathbf{b}(K)$ has such a gap, the corresponding $(g,b)$-surface is irreducible. \\

It is well known that the bridge spectra of torus knots have a single gap.  In the present work, we set out to answer the following question, first proposed to us by Yo'av Rieck:

\begin{question}\label{quest}
Are there examples of knots in $S^3$ whose bridge spectra have more than one gap?
\end{question}

A theorem of Tomova completely characterizes the bridge spectra of high distance knots:

\begin{theorem} \cite{tomova}
Suppose $K$ is a knot in $S^3$ with a $(0,b)$-bridge sphere $\Sigma$ of sufficiently high distance (with respect to $b$).  Then any $(g',b')$-bridge surface $\Sigma'$ satisfying $b' = b_{g'}(K)$ is the result of meridional stabilizations performed on $\Sigma$.  Thus
\[ \mathbf{b}(K) = (b,b-1,b-2,\dots,0).\]
\end{theorem}

It follows that the bridge spectrum of a ``generic" knot $K$ is rather uninteresting.  We will show that, in contrast to high distance knots, the bridge spectra of iterated torus knots exhibit different behavior.  The main theorem is as follows:

 \begin{theorem}\label{main}
 Suppose that $K_n = ((p_0,q_0),\dots,(p_n,q_n))$ is an iterated torus knot, whose cabling parameters satisfy $|p_i - p_{i-1}q_{i-1}q_i| > 1$.  Then
 \[b_g(K_n)=	\begin{cases}
			q_n \cdot b_g(K_{n-1}) &\text{if $g < n$;}\\
			\min\{|p_n - p_{n-1}q_{n-1}q_n|,q_n\} &\text{if $g = n$;}\\
			0 &\text{otherwise.}
		\end{cases}
\]
In other words,
\[ \mathbf{b}(K_n) = q_n \cdot \mathbf{b}(K_{n-1}) + \min\{|p_n - p_{n-1}q_{n-1}q_n|,q_n\} \cdot \mathbf{e}_n.\]
\end{theorem}

It follows that the bridge spectrum of $K_n$ has a gap at every index from 1 to $n+1$, providing a positive answer to Question \ref{quest} above and yielding for any $n$ the first examples of a knot $K_n$ in $S^3$ having more than $n$ irreducible bridge surfaces.  In the course of proving the main theorem, we show another result, previously unknown for $n>1$: the tunnel number of such $K_n$ is $n+1$.  This proof uses a theorem of Schultens concerning Heegaard splittings of graph manifolds \cite{schult2}. \\

The main theorem is related to a classical result of Schubert \cite{schub}, with a modern proof given by Schultens \cite{schult}:

\begin{theorem}\cite{schub}\cite{schult}\label{satbridge}
Let $K$ be a satellite knot with companion $J$ and pattern of index $n$.  Then
\[ b_0(K) \geq n \cdot b_0(J).\]
\end{theorem}
A $(g,b)$-bridge surface for a knot $K$ is \emph{minimal} if $b = b_g(K)$.  Schultens' proof of the theorem reveals that after isotopy, a minimal bridge sphere $\Sigma$ for $K$ can be made to intersect the companion torus in meridian disks, so that each bridge of $J$ contributes at least $n$ bridges to $K$.  The proof of Theorem \ref{main} yields that whenever $g < n$, a minimal genus $g$ bridge surface for $K_n$ can be made to intersect the companion torus corresponding to $K_{n-1}$ in meridian disks.  However, there are knots $K_n$ such that a minimal genus $n$ bridge surface meets the companion torus in an annulus; hence, we cannot, in general, hope for an extension of Theorem \ref{satbridge} to surfaces of higher genus. \\

The proof of Theorem \ref{main} requires a lemma regarding strongly irreducible bridge surfaces, which may be of independent interest.  Roughly, this lemma says that a strongly irreducible bridge surface cut along properly embedded essential surfaces will decompose into at most one exceptional component which is strongly irreducible, along with some number of incompressible components.

\begin{lemma}
Let $M$ be a compact, orientable, irreducible 3-manifold, $J \subset M$ a properly embedded 1-manifold, and $Q = \pd N(J) \subset \pd M(J)$.  Suppose $\Sigma$ is a strongly irreducible bridge splitting surface for $(M,J)$, and let $S \subset M(J)$ be a collection of properly embedded essential surfaces such that for each component $c$ of $\pd S$, either $c \subset Q$ or $c \subset \pd M$.  Then one of the following must hold:

\be
\item[(1)] After isotopy, $\Sigma_J$ is transverse to $S$ and each component of $\Sigma_J \setminus \eta(S)$ is $Q$-essential in $M(J) \setminus \eta(S)$.

\item[(2)] After isotopy, $\Sigma_J$ is transverse to $S$, one component of $\Sigma_J \setminus \eta(S)$ is $Q$-strongly irreducible and all other components are $Q$-essential in $M(J) \setminus \eta(S)$,

\item[(3)] After isotopy, $\Sigma_J$ is almost transverse to $S$, and each component of $\Sigma_J \setminus \eta(S)$ is $Q$-essential in $M(J) \setminus \eta(S)$.

\ee
\end{lemma}

Combined with results on strongly irreducible bridge surfaces by Hayashi and Shimokawa \cite{hayshi}, the above lemma can easily be seen to provide alternate proofs of the theorems of \cite{ozawa} and \cite{zupan} regarding the destabilizability of torus knots and cables of $mp$-small destabilizable knots. \\

The article is organized as follows:  In Section \ref{pre}, we introduce relevant definitions and background material.  Section \ref{iter} discusses topological properties of iterated torus knots, while Section \ref{tunnel} contains a complete analysis of the tunnel number of such knots.  Section \ref{strong} presents the above lemma about strongly irreducible bridge surfaces, setting up the proof of the main theorem in Section \ref{spec}.  Section \ref{cableex} includes a worked example with figures, and finally Section \ref{question} poses some open problems that may be of interest. \\

\noindent \textbf{Acknowledgements} \, I would like to thank Marion Campisi, Cameron Gordon, Tye Lidman, and Maggy Tomova for helpful conversations and insights over the course of this project.

\section{Preliminaries}\label{pre}

Throughout, all 3-manifolds and surfaces will be compact, orientable, and irreducible.  We will let $\eta( \cdot)$ and $N( \cdot)$ denote open and closed regular neighborhoods, respectively, in an ambient manifold that should be clear from context.  Let $S$ be a properly embedded surface in a 3-manifold $M$.  A \emph{compressing disk} $D$ for $S$ is an embedded disk such that $D \cap S = \pd C$ but $\pd D$ does not bound a disk $D' \subset S$.  A \emph{$\pd$-compressing disk} $\Delta$ for $S$ is an embedded disk such that $\pd \Delta$ is the endpoint union of arcs $\gamma_1$ and $\gamma_2$ such that $\Delta \cap S = \gamma_1$, $\gamma_1$ is essential in $S$, and $\Delta \cap \pd M = \gamma_2$.  The surface $S$ is said to be \emph{incompressible} if there does not exist a compressing disk $D$ for $S$ and \emph{$\pd$-incompressible} if there does not exist a $\pd$-compressing disk $\Delta$ for $S$.  Further, $S$ is \emph{essential} if $S$ is incompressible, $\pd$-incompressible, and not parallel into $\pd M$. \\

Suppose now that $M$ a 3-manifold containing a properly embedded 1-manifold $J$, and denote the \emph{exterior} of $J$ in $M$ by $M(J) = M \setminus \eta(J)$ (if $M = S^3$, we write $E(J) = M(J)$).  Let $\Sigma$ be a properly embedded surface transverse to $J$, with $\Sigma_J$ denoting $\Sigma \setminus \eta(J)$.  A \emph{compressing disk} $D$ for $\Sigma_J$ in $(M,J)$ is an embedded disk $D$ such that $D \cap \Sigma_J = \pd D$, $D \cap J = \emp$, and $\pd D$ does not bound a disk in $\Sigma_J$.  A \emph{bridge disk} $\Delta$ for $\Sigma_J$ in $(M,J)$ is an embedded disk such that $\pd \Delta$ is the endpoint union of two essential arcs $\gamma_1$ and $\gamma_2$, where $\Delta \cap \Sigma_J = \gamma_1$ and $\Delta \cap J = \gamma_2$.  We also think of $\Delta$ as a $\pd$-compressing disk for $\Sigma_J$ in $M(J)$ with $\Delta \cap \pd M(J) \subset \pd N(J)$.  Finally, a \emph{cut disk} $E$ for $\Sigma_J$ in $(M,J)$ is an embedded disk $C$ such that $C \cap \Sigma_J = \pd C$, $C \cap J$ is a single point, and $\pd C$ does not bound a disk in $\Sigma_J$. \\

Let $V$ be a compression body and $\A \subset V$ a collection of properly embedded arcs.  We say that $\A$ is \emph{trivial} if every arc is either vertical or isotopic into $\pd_+V$.  A \emph{bridge splitting} of $(M,J)$ with \emph{bridge surface} $\Sigma$ is the decomposition of $(M,J)$ as $(V,\A) \cup_{\Sigma} (W,\n)$, where each $V$ and $W$ are compression bodies containing collections $\A$ and $\n$ of trivial arcs.  In the special case that $\A$ and $\n$ contain only boundary parallel arcs, we say that $\Sigma$ is a $(g,b)$\emph{-bridge splitting surface}, where $g = g(\Sigma)$ and $b=|\A| = |\n|$. \\

Given a bridge splitting $(M,J) = (V,\A) \cup_{\Sigma} (W,\n)$, there are several generic methods to construct new bridge surfaces for $(M,J)$.  To increase the genus of the splitting, let $\gamma \subset V$ be a $\pd_+$-parallel arc such that $\gamma \cap \A = \emp$.  Defining $V' = V \setminus \eta(\gamma)$, $W' = W \cup N(\gamma)$, and $\Sigma' = \pd_+ V' = \pd_+ W'$, we have $(M,J) = (V',\A) \cup_{\Sigma'} (W',\n)$ is another splitting of $(M,J)$.  This process is called \emph{elementary stabilization}.  Note that the surface $\Sigma'$ has compressing disks that intersect in a single point.  If, in the reverse direction, there are compressing disks $D$ in $(V,\A$) and $E$ in $(W,\n)$ for $\Sigma_J$ such that $|D \cap E| =1$, then $\pd N(D \cup E)$ is a 2-sphere which bounds a ball and which intersects $\Sigma$ in a single curve.  Compressing $\Sigma$ along $\Sigma \cap \pd (N(D \cup E))$ yields a bridge surface $\Sigma''$ of lower genus, and we say that $\Sigma$ is \emph{stabilized}. \\

We may also increase the number of trivial arcs in the splitting by adding an extra pair of canceling trivial arcs to $\A$ and $\n$ near some point of $J \cap \Sigma$.  The resulting surface $\Sigma'$ is called an \emph{elementary perturbation} of $\Sigma$.  Conversely, if there are bridge disks $\Delta$ and $\Delta'$ in $(V,\A)$ and $(W,\n)$ such that $\Delta \cap \Delta'$ is a single point contained in $J$, we may construct an isotopy which cancels two arcs of $\A$ and $\n$, creating a new surface $\Sigma''$, and we say that $\Sigma$ is \emph{perturbed}.  If $\Sigma$ is a $(g,b)$-splitting surface and $\Sigma_1$ is an elementary stabilization of $\Sigma$, then $\Sigma_1$ is a $(g+1,b)$-surface.  If $\Sigma_2$ is an elementary perturbation of $\Sigma$, then $\Sigma_2$ is a $(g,b+1)$-surface.  Given any two bridge surfaces $\Sigma_1$ and $\Sigma_2$ for $(M,J)$, there is a third bridge surface $\Sigma^*$ which can be obtained from either $\Sigma_i$ by elementary perturbations and stabilizations \cite{zupan1}. \\

If there are bridge disks $\Delta$ and $\Delta'$ in $(V,\A)$ and $(W,\n)$ such that $\Delta \cap \Delta'$ is two points contained in $J$, then a component of $J$ is isotopic into $\Sigma$ and we say $\Sigma$ is \emph{cancelable}.  Here we set the convention that a $(g,0)$-bridge surface for a knot $K$ in $M$ is a Heegaard surface $\Sigma$ for $M$ such that $K \subset \Sigma$.  In some settings, it is required in addition that $K$ be primitive with respect to one of the handlebodies, but we do not make that restriction here.  Note that a $(g,0)$-surface may be perturbed to a $(g,1)$-surface. \\

A slightly more complicated way to construct a new surface from $\Sigma$ is the following: fix a $\pd_+$-parallel component $\A_1$ of $\A$.  Letting $V' = V \setminus \eta(\A_1)$ and $W' = W \cup N(\A_1)$, we have that $\Sigma' = \pd V' = \pd W'$ is a bridge surface for $(M,J)$.  We call this process \emph{meridional stabilization}.  If $\Sigma$ is a $(g,b)$-surface, the resulting surface $\Sigma'$ will be a $(g+1,b-1)$-suface.  On the other hand, if there exists a compressing disk $D$ for $\Sigma_J$ in $(V,\A)$ and a cut disk $C$ for $\Sigma_J$ in $(W,\n)$ (or vice versa) such that $|C \cap D| = 1$, we say that $\Sigma$ is \emph{meridionally stabilized} and we may reverse the above process to construct a new splitting $\Sigma''$ of $(M,J)$. \\

Recall that the \emph{genus $g$ bridge number} $b_g(K)$ of a knot $K \subset S_3$ is defined as
\[ b_g(K) = \min\{b : K \text{ admits a } (g,b)\text{-splitting}\}.\]
Observe that $b_0(K)$ is the classical bridge number of $K$.  Further, using meridional stabilization, we have $b_{g+1}(K) \leq b_g(K) - 1$.  It follows that $b_g(K) = 0$ for $g \geq b_0(K)$. \\ 

In order to consider all genus $g$ bridge numbers of $K$ simultaneously, we define the \emph{bridge spectrum} $\mathbf{b}(K)$ of $K$:
\[ \mathbf{b}(K) = (b_0(K),b_1(K),b_2(K),\dots)\]
Although $\mathbf{b}(K) \in \Z^{\infty}$, the above argument implies that $\mathbf{b}(K)$ must have finitely many nonzero entries, and $\mathbf{b}(K)$ is bounded above by the sequence $(b_0(K),b_0(K)-1,b_0(K)-2,\dots)$. \\

A $(g,b)$-surface for a knot $K$ is said to be \emph{irreducible} if it is not stabilized, perturbed, meridionally stabilized, or cancelable.  From the point of view of the bridge spectrum, if $b_{g}(K) < b_{g-1}(K) - 1$, a genus $(g,b)$-surface $\Sigma$ satisfying $b = b_g(K)$ must be irreducible.  For this reason, suppose $b_{g}(K) < b_{g-1}(K)- 1$ and let $n = b_{g-1}(K) - b_{g}(K)$.  In this case, we say that $\mathbf{b}(K)$ has a \emph{gap of order $n$ at index $g$}. \\

%In the following theorems, a Heegaard splitting is sufficiently complicated if it has Hempel distance above a certain threshold related to the genus, and a $(g,b)$-bridge splitting $(S^3,K) = (V,\A) \cup_{\Sigma} (W,\n)$ is \emph{sufficiently complicated} if the distance in the curve complex between the disk sets of $(V,\A)$ and $(W,\n)$ is above a threshold determined by $g$ and $b$. \\

Several results over $(g,b)$-splittings can be adapted to statements about bridge spectra.  In \cite{mms}, Minsky, Moriah, and Schleimer prove a theorem which we can state as follows:

\begin{theorem}\label{MMS}
For every $n$ and $m$ such that $n \geq 2$ and $m \geq 1$, there exist knots $K$ such that $\mathbf{b}(K)$ has a gap of order at least $n$ at index $m$.
\end{theorem}

%there exist knots $K \subset S^3$ with the property that $b_{m}(K) = 0$, while $\mathbf{b}(K)$ is bounded below by $(n+m-1,n+m-2,\dots,n,0)$.

The knots produced in Theorem \ref{MMS} have the property that $b_m(K) = 0$.  More generally, in \cite{tomova} Tomova shows that
\begin{theorem}\label{tom}
For every $n$, $m$, and $l$, where $n \geq 2$ and $m,l \geq 1$, there exist knots $K$ such that $\mathbf{b}(K)$ has a gap of order at least $n$ at index $m$ and $b_m(K) = l$.  Further, there exist knots $K$ whose bridge spectra have no gaps.
\end{theorem}

%there exist knots $K \subset S^3$ with the property that $b_g(K) = l+m-g$ whenever $g \geq m$, while $\mathbf{b}(K)$ is bounded below by $(n+m+l-1,n+m+l-2,\dots,n+l,l,l-1,\dots,0)$.  In addition, for $n=1$ and $m=0$, there are knots $K$ such that $\mathbf{b}(K) = (l,l-1,\dots,2,1,0)$.

The knots $K$ from Theorem \ref{tom} whose spectra have gaps at index $m$ have the additional property that $b_g(K) = l+ m - g$ for $g \geq m$, and corresponding minimal surfaces are meridional stabilizations of an $(m,l)$-surface.  We note that the above examples are produced by exhibiting sufficiently complicated bridge and Heegaard surfaces, where complexity is measured by examining distance between disks sets in the curve complex of the surface.  These methods, however, are not suited to answer Question \ref{quest} regarding knots whose spectra have more than one gap. \\

In order to classify bridge spectra, we employ generalizations of Heegaard splittings developed by Scharlemann and Thompson \cite{scharthom} and adapted to bridge splittings by Hayashi and Shimokawa \cite{hayshi}.  Let $M$ be a 3-manifold containing a properly embedded 1-manifold $J$.  We will describe the theories of generalized Heegaard and bridge splittings simultaneously by using the bridge splitting terminology and noting all results hold in the case that $J = \emp$, in which case bridge surfaces become Heegaard surfaces. \\

Let $\Sigma$ be a bridge surface for $(M,J)$ which yields the splitting $(M,J) = (V,\A) \cup_{\Sigma} (W,\n)$.  We say that $\Sigma$ is \emph{weakly reducible} if there exist disjoint compressing or bridge disks $D \subset (V,\A)$ and $D' \subset (W,\n)$ for $\Sigma_J$.  If $\Sigma$ is not weakly reducible, perturbed, or cancelable, we say $\Sigma$ is \emph{strongly irreducible}.  By considering bridge disks as embedded in $M(J)$, we can see that perturbed and cancelable surfaces will be weakly reducible; hence, in $M(J)$, $\Sigma$ is strongly irreducible if and only if it is not weakly reducible. \\

Now, suppose $(M,J)$ contains a collection $\mathcal{S}$ $=$ $\{\Sigma_0,S_1,\Sigma_1,\dots,S_d,\Sigma_d\}$ of disjoint surfaces transverse to $J$ and such that $(M,J)$ cut along $\mathcal{S}$ is a collection of compression bodies containing trivial arcs $\{(C_0,\tau_0),(C_0',\tau'_0),\dots,$ $(C_d,\tau_d),(C_d',\tau'_d)\}$, where 
\bi
\item $(C_i,\tau_i) \cup_{\Sigma_i} (C_i',\tau_i')$ is a bridge splitting of some submanifold $(M_i,J_i)$ for each $i$,
\item $\pd_- C_i = \pd_- C_{i-1}' = S_i$ for $1 \leq i \leq d$,
\item $\pd M = \pd_- C_0 \cup \pd_- C_d'$, and
\item $J = \bigcup (\tau_i \cup \tau_i')$
\ei

We call this decomposition a \emph{multiple bridge splitting}, and the surfaces $\Sigma_i$ are called \emph{thick} whereas the surfaces $S_j$ are called $\emph{thin}$.  The thick surface $\Sigma_i$ is strongly irreducible if it is strongly irreducible in the manifold $(C_i,\tau_i) \cup_{\Sigma_i} (C_i',\tau_i')$.  A multiple splitting is called \emph{strongly irreducible} if each thick surface is strongly irreducible and no compression body is trivial (homeomorphic to $\Sigma_i \X I$ with $\tau_i$ only vertical arcs).  The following crucial theorem comes from Hayashi and Shimokawa \cite{hayshi} (and in the special case that $J = \emp$, it is proved by Scharlemann and Thompson \cite{scharthom}).  Other proofs of this fact are given in the contexts of $\A$-sloped Heegaard surfaces \cite{campisi} and in the more general setting of embedded graphs \cite{taylortomova}.

\begin{theorem}\label{thin}
Let $M$ be a 3-manifold containing a 1-manifold $J$.  If $(M,J)$ has a strongly irreducible multiple bridge splitting, then $\pd M_J$ and every thin surface is incompressible.  On the other hand, if $\pd M_J$ is incompressible in $M(J)$ and $\Sigma$ is a weakly reducible bridge splitting for $(M,J)$, then $(M,J)$ has a strongly irreducible multiple bridge splitting satisfying
\begin{equation}\label{amal}
g(\Sigma) = \sum g(\Sigma_i) - \sum g(S_i).
\end{equation}
\end{theorem}
Although we will not go into details here, the multiple splitting given by the above theorem is obtained from the weakly reducible surface $\Sigma$ via a process known as \emph{untelescoping}.  Conversely, if $\mathcal{S}$ yields a multiple bridge splitting, we may construct a bridge surface $\Sigma$ satisfying (\ref{amal}) via a process called \emph{amalgamation}.

\section{Iterated torus knots}\label{iter}

We will focus on a class of knots known as iterated torus knots, which make up a subset of a collection of knots called cable knots.  Let $V$ be a unknotted torus standardly embedded in $S^3$, with $T = \pd V$.  For a curve $c \subset T$, $[c]$ may be parametrized as $(p,q) = p[\mu] + q[\lambda]$ in $H_1(T)$, where $\mu$ bounds a meridian disk of $V$ and $\lambda$ is the preferred longitude of $V$, bounding a disk outside of $V$.  Choosing $q \geq 2$, let $\hat{K}_{p,q}$ be a copy of a $(p,q)$-curve on $T$ pushed into $\text{int}(V)$.  Now, suppose that $K$ is a knot in $S^3$, and let $\varphi(V) \rightarrow S^3$ be a knotted embedding of $V$ such that if $C$ is a core of $V$, then $\varphi(C)$ is isotopic to $K$.  In this case $K_{p,q} = \varphi(\hat{K}_{p,q})$ is called a \emph{$(p,q)$-cable} of $K$.   In addition, we stipulate that the cable has the preferred framing; that is, we have $[\varphi(\lambda)]=0$ in $H_1(E(K))$. \\

\begin{figure}[h!]
  \centering
    \includegraphics[width=.6\textwidth]{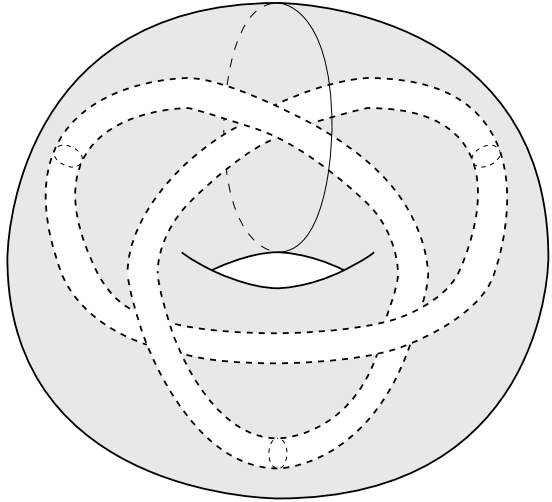}
    \caption{The cable space $C_{3,2}$}
    \label{cablpac}
\end{figure}

Using the same symbols as above, we follow \cite{lither} to define a $(p,q)$\emph{-cable space} $C_{p,q}$ by
\[ C_{p,q} = V \setminus \eta(\hat{K}_{p,q}).\]
Observe that $E(K_{p,q})$ decomposes as $E(K) \cup_T C_{p,q}$, where we consider the torus $T$ as both $\pd V$ and $\pd E(K)$.  The space $C_{p,q}$ is Seifert fibered with base space an annulus and a single exceptional fiber $C$ (for further discussion, see \cite{scott}).  We distinguish between the outer boundary $\pd_+ C_{p,q} = \pd V$ and the inner boundary $\pd_- C_{p,q} =  \pd N(\hat{K}_{p,q})$.  Note that $H_1(\pd_- C_{p,q})$ inherits a natural basis by virtue of $\hat{K}_{p,q}$ being a torus knot.  See Figure \ref{cablpac}. \\

The class of \emph{iterated torus knots} is defined inductively:  Define $K_0$ to be a $(p_0,q_0)$-torus knot.  Then for any $n \geq 1$, define $K_n$ to be a $(p_n,q_n)$-cable of $K_{n-1}$.  We call the knot $K_n$ an \emph{iterated torus knot}, denoted by $K_n = ((p_0,q_0),\dots,(p_n,q_n))$.  Here $E(K_n)$ decomposes as
\[ E(K_n) = E(K_0) \cup_{T_1} C_{p_1,q_1} \cup_{T_2} \dots \cup_{T_n} C_{p_n,q_n},\]
where $T_i = \pd_+ C_i = \pd_- C_{i-1}$ for $i \geq 2$ and $T_1 = \pd_+ C_1 = \pd E(K_0)$.  The longitude-meridian bases of $H_1(\pd_+ C_i)$ and $H_1(\pd_- C_{i-1})$ coincide; hence, we specify a preferred basis $(\mu_i,\lambda_i)$ for $H_1(T_i)$ without ambiguity.  For a curve $c \subset T_i$ with $[c] = a[\mu_i] + b[\lambda_i]$, we say that $c$ has \emph{slope} $\frac{a}{b}$. \\

By \cite{hatcher}, if $M$ is a Seifert fibered space and $S$ is incompressible in $S$, then $S$ must either be a union of fibers, called a \emph{vertical surface}, or transverse to all fibers, called a \emph{horizontal surface}.  There is a natural projection map $\pi: M \rightarrow B$, where $B$ is the base orbifold of $M$, and if $S$ is a horizontal surface, $\pi |_S$ is an orbifold covering map. \\

The following lemma, due to Gordon and Litherland \cite{lither}, provides a complete classification of incompressible surfaces in cable spaces.  It will be critical for our understanding of such surfaces in $E(K_n)$.  For two curves with slopes $\frac{a}{b}$ and $\frac{c}{d}$, we will denote their intersection number by $\Delta(\frac{a}{b},\frac{c}{d}) = |ad-bc|$.

\begin{lemma}\label{L1}
Suppose that $S \subset C_{p,q}$ is incompressible and not $\pd$-parallel.
\be
\item If $S$ is vertical, then $S$ is an annulus, and any component of $S \cap \pd_+ C_{p,q}$ has slope $\frac{p}{q}$, while any component of $S \cap \pd_- C_{p,q}$ has slope $pq$.
\item If $S$ is horizontal, then there are coprime integers $m$ and $n$ (with $n\geq 0$) such that the total boundary slope of $S \cap \pd_+ C_{p,q}$ is $\frac{n+mp}{qm}$ and the total boundary slope of $S \cap \pd_- C_{p,q}$ is $\frac{q(n+mp)}{m}$.  In this case, $\chi(S) = n(1-q)$.
\item If $S$ is horizontal and planar, then either
{\be 
\item there are $q$ components of $S \cap \pd_+ C_{p,q}$, each of slope $\frac{l}{m}$, and one component of $S \cap \pd_- C_{p,q}$, of slope $\frac{lq^2}{m}$, where $l$ and $m$ satisfy $\Delta(\frac{l}{m},\frac{p}{q}) = 1$; or
\item there is one component of $S \cap \pd_+ C_{p,q}$, of slope $\frac{1+kpq}{kq^2}$, and $q$ components of $S \cap \pd_- C_{p,q}$, each of slope $\frac{1 + kpq}{k}$, where $k$ is an integer.
\ee}
\ee
\end{lemma}

We may employ Lemma \ref{L1} in the proof of the statements below.  Note that for any cable knot $K_{p,q}$, the companion torus $T$ is essential in $E(K_{p,q}$); hence in the iterated torus knot $K_n$, each torus $T_i$ is an essential surface.  

\begin{lemma}\label{L2}
Suppose $S \subset C_{p,q}$ is incompressible.  If each component of $S \cap \pd_+ C_{p,q}$ has integral slope, then each component of $S \cap \pd_- C_{p,q}$ also has integral slope.
\begin{proof}
Suppose without loss of generality that $S$ is connected.  Since $S \cap \pd_+ C_{p,q}$ has integral slope, $S$ is not vertical, and so by Lemma \ref{L1}, the total slope of $S \cap \pd_+ C_{p,q}$ is $\frac{n+mp}{qm}$ for some coprime integers $m$ and $n$; that is, there are $r = \text{gcd}(n+mp,qm)$ components of $S \cap \pd_+ C_{p,q}$, each of slope $\frac{a}{b}$, where $a = \frac{n+mp}{r}$ and $b = \frac{qm}{r}$.  By assumption, $b = \pm 1$; thus $r = \pm qm$ and $m$ divides $n+ mp$.  It follows that $m$ divides $n$; hence $m = \pm 1$.  The total slope of $S \cap \pd_- C_{p,q}$ is $\frac{q(m+np)}{m}$ and as such is also integral.
\end{proof}
\end{lemma}

A boundary component of a surface is said to be \emph{meridional} if its boundary slope is $\frac{1}{0}$.

\begin{lemma}\label{L3}
Suppose $S \subset C_{p,q}$ is incompressible.  Then $S \cap \pd_+ C_{p,q}$ is meridional if and only if $S \cap \pd_- C_{p,q}$ is meridional.
\begin{proof}
It is clear that if either set of boundary components is meridional, then $S$ is horizontal.  In this case, by Lemma \ref{L1}, the total outer boundary slope is $\frac{n+mp}{qm}$ while the total inner boundary slope is $\frac{q(n+mp)}{m}$.  If $m=0$, then both outer and inner boundary components of $S$ are meridional; otherwise, neither is meridional.
\end{proof}
\end{lemma}

In light of Lemma \ref{L3}, we say that a horizontal surface $S$ in $C_{p,q}$ is \emph{meridional} if any component of $\pd S$ is meridional.

\begin{lemma}\label{L4}
Suppose $S \subset C_{p,q}$ is incompressible.  If each component of $S \cap \pd_+ C_{p,q}$ is a single component of slope $\frac{1}{n}$ and $S$ is planar, then $S$ is meridional.
\begin{proof}
By Lemma \ref{L1}, we have $\frac{1}{n} = \frac{1 + kpq}{kq^2}$; hence $k=0$.
\end{proof}
\end{lemma}

Let $\mathcal{T}$ denote the union of the tori $T_i$ and let $T_{n+1}$ denote $\pd E(K_n)$.

\begin{lemma}\label{L5}
Suppose $S \subset E(K_n)$ is an essential surface.  If $S$ is not isotopic to some $T_i$, then $S \cap \pd E(K_n)$ is nonempty and has integral slope.
\begin{proof}
After isotopy, we may suppose that $|S \cap \mathcal{T}|$ is minimal.  First, we note that each component of $S \cap C_i$ or $S \cap E(K_0)$ is incompressible in $C_i$ or $E(K_0)$:  If $D$ is a compressing disk for $S \cap C_i$ in $C_i$, then $\pd D$ bounds a disk $D' \subset S$ by the incompressibility of $S$, where $D' \cap \mathcal{T} \neq \emp$.  By the irreducibility of $E(K_n)$, there is an isotopy of $S$ pushing $D'$ onto $D$ which reduces $|S \cap \mathcal{T}|$, a contradiction. \\

If $S \cap \mathcal{T} = \emp$, then $S \subset C_i$ for some $i$ or $S \subset E(K_0)$.  However, in this case $S$ must be vertical and closed, and the only such surfaces in $C_i$ or $E(K_0)$ are boundary parallel.  It follows that $S$ is isotopic to $T_i$ for some $i$. \\

Thus, suppose that $S \cap \mathcal{T} \neq \emp$.  If $S \cap T_1 \neq \emp$, then $S \cap E(K_0)$ is either a Seifert surface with boundary slope $0$ or an essential annulus with slope $p_0q_0$.  If $S \cap T_1 = \emp$, let $i$ be the smallest $i$ such that $S \cap T_i \neq \emp$.  Lemma \ref{L1} dictates that $S \cap C_{i-1}$ must be a vertical annulus in $C_{i-1}$ which intersects $T_i$ in a curve with slope $p_iq_i$ (otherwise $S \cap T_{i-1} \neq \emp$).  In either case we have that $S \cap T_i$ has integral slope.  By a repeated application of Lemma \ref{L2}, $S \cap T_j$ has integral slope for all $j \geq i$, completing the proof.

\end{proof}
\end{lemma}

\section{The tunnel number of iterated torus knots}\label{tunnel}

The \emph{tunnel number} $t(K)$ of a knot $K \subset S^3$ is defined as $t(K) = g(E(K)) - 1$, where $g(E(K))$ is the Heegaard genus of $E(K)$.  Our goal in this section is to determine the minimal genus of a Heegaard splitting of $E(K_n)$ for an iterated torus knot $K_n$, and our main tool is the classification of a Heegaard splittings of graph manifolds by Schultens \cite{schult2}. \\

A \emph{graph manifold} is a 3-manifold $M$ containing a nonempty collection of essential tori $T$ such that $M \setminus \eta(T)$ is a disjoint union of Seifert fibered spaces.  A surface $S$ properly embedded in a Seifert fibered space $M'$ is \emph{pseudohorizontal} if after isotopy there exists a fiber $f$ of $M'$ such that $S \cap (M' \setminus \eta(f))$ is a horizontal surface in the space $M' \setminus \eta(f)$, and $S \cap N(f)$ is an annulus which is a bicollar of $f$.  Note that each boundary component of a horizontal surface has the same orientation; thus, if a pseudohorizontal surface $S$ is to be orientable, then the horizontal piece $S \cap (M' \setminus \eta(f))$ must have two components. \\

Suppose that $S$ and $S'$ are connected surfaces properly embedded in $M$ and $\gamma$ is an arc with one endpoint in each of $S$ and $S'$ and interior disjoint from $S \cup S'$.  Then $N(S \cup S' \cup \gamma)$ has three boundary components, one isotopic to $S$, another isotopic to $S'$, and third we define to be the result of \emph{tubing $S$ to $S'$ along $\gamma$}. \\

We may now state one of Schultens' main results over the classification of graph manifolds \cite{schult2}:

\begin{theorem}\label{graph}
Suppose $M$ is a graph manifold containing a strongly irreducible Heegaard surface $\Sigma$, and let $T$ be a collection of essential tori splitting $M$ into Seifert fibered spaces.  After an isotopy of $\Sigma$, there is a submanifold $M'$ of $T$ such either $M'$ is a component of $M \setminus \eta(T)$ or $M = T' \X I$ for a torus $T' \in T$, called the active component of $M$, with the following properties:
{\be
\item Each component of $\Sigma \setminus (\eta(T) \cup \eta(M'))$ is incompressible in $M \setminus (\eta(T) \cup \eta(M'))$.
\item If $M' \neq T' \X I$ for some $T' \in T$, then $\Sigma \cap M'$ is a pseudohorizontal surface in $M'$.
\item Otherwise, $M' = T' \X I$, and $\Sigma \cap M'$ is either
{\be
\item the result of tubing an annulus parallel into $T' \X \{0\}$ with boundary slope $s_0$ to an annulus parallel into $T' \X \{1\}$ with boundary slope $s_1$ along a vertical arc in $T' \X I$, where $\Delta(s_0,s_1) = 1$, or
\item the result of tubing two vertical annuli in $T' \X I$ along an arc contained in $T \X \{\frac{1}{2}\}$.
\ee}
\ee}
\end{theorem}

From this point forward, for an iterated torus knot $K_n = ((p_0,q_0),(p_1,q_1),$ $\dots,(p_n,q_n))$, we make the restriction that $|p_i - p_{i-1}q_{i-1}q_i| > 1$ for all $i$.  The motivation for this is as follows:  A regular fiber $f_{i-1}$ of $C_{i-1}$ has slope $s_{i-1} = p_{i-1}q_{i-1}$ in $T_i$, whereas a regular fiber $f_i$ of $C_i$ has slope $s_i = \frac{p_i}{q_i}$ in $T_i$.  In order for the arguments below to hold, we need that $|f_i \cap f_{i-1}| > 1$.  In other words, we require
\[ |p_i - p_{i-1}q_{i-1}q_i| = \Delta(s_i,s_{i-1}) > 1.\]

The goal of this section is to show that for such iterated torus knots $K_n$, we have $t(K_n) = n+1$.  We will use Theorem \ref{thin} in conjunction with Theorem \ref{graph} to bound the genus of any Heegaard surface for $E(K_n)$ from below.  For this, we must understand pseudohorizontal surfaces in $E(K_0)$ and the cable spaces $C_i$.

\begin{lemma}\label{L6}
Suppose that $S$ is a pseudohorizontal surface in $E(K_0)$.  Then either $\chi(S) \leq -4$ or $\pd S$ has meridional slope.
\begin{proof}
Without loss of generality, let $p_0 \geq 2$ and let $f$ be a fiber of $E(K_0)$ such that $S' = S \cap (E(K_0) \setminus \eta(f))$ is horizontal.  First, suppose that $f$ is a regular fiber, and let $B'$ be the base space of the Seifert fibered space $E(K_0) \setminus \eta(f)$, so that $B'$ is a topological annulus with two cone points.  By the observation above, $S'$ has two components $S_1$ and $S_2$,  each an orbifold cover of $B'$ of degree $d \geq p_0q_0$.  Following the discussion in \cite{scott}, for instance, we let $B^*$ denote $B'$ with a neighborhood of the cone points removed; thus
\begin{eqnarray*}
\chi(S_i) &=& d \chi(B^*) + \frac{d}{p_0} + \frac{d}{q_0} \\
&\leq& -2 p_0q_0 + p_0 + q_0 \\
&=&  p_0(1-q_0) + q_0(1-p_0) \leq -7.
\end{eqnarray*}
It follows that $\chi(S) = \chi(S') \leq -14$. \\

In the second case, suppose that $f$ is a critical fiber.  Without loss of generality, we may suppose that $E(K_0) \setminus \eta(f) = C_{p_0,q_0}$.  After isotopy, $S$ intersects $N(f)$ in an annulus which is a bicollar of $f$; thus $S \cap N(f)$ has integral boundary slope on the solid torus $N(f)$.  In the coordinates of $\pd_+ C_{p_0,q_0}$, each component of $S \cap \pd_+ C_{p_0,q_0}$ has slope $\frac{1}{k}$ for some integer $k$.  Further, $S'$ has two components $S_1$ and $S_2$, each with a single boundary component in $\pd_+ C_{p_0,q_0}$.  By Lemma \ref{L4}, if each $S_i$ is planar, then $S$ has meridional boundary.  Otherwise, $\chi(S_i) \leq -2$, and $\chi(S) = \chi(S') \leq -4$, as desired.
\end{proof}
\end{lemma}

We carry out a similar analysis of pseudohorizontal surfaces in a cable space $C_{p,q}$:

\begin{lemma}\label{L7}
Suppose that $S$ is a pseudohorizontal surface in the cable space $C_{p,q}$.  Then $\chi(S) \leq -4$.  If, in addition, each component of $S \cap \pd_+ C_{p,q}$ has boundary slope $\frac{l}{1}$, where $\Delta(l,\frac{p}{q}) > 1$, then either $\chi(S) \leq -6$, the slope of $S \cap \pd_- C_{p,q}$ is meridional, or the slope of $S \cap \pd_- C_{p,q}$ is integral.
\begin{proof}
Let $f$ be a fiber of $C_{p,q}$ such that $S' = S \cap (C_{p,q} \setminus \eta(f))$ is horizontal.  As above, we first suppose that $f$ is a regular fiber, and let $B'$ be the base space of $C_{p,q} \setminus \eta(f)$, where $B'$ is a thrice-punctured sphere with one exceptional fiber.  Again, $S'$ has two components $S_1$ and $S_2$, each of which is a degree $d \geq q$ orbifold cover of $B'$, so that
\[ \chi(S_i) = -2d + \frac{d}{q} = d\left(-2 + \frac{1}{q}\right) \leq -2q  +1 \leq -3,\]
so $\chi(S) = \chi(S') \leq -6$. \\

Otherwise, suppose $f$ is the unique exceptional fiber of $C_{p,q}$, and note that $H_1(\pd N(f))$ inherits the same natural basis as $H_1(\pd_+ C_{p,q})$.  Letting $S' = S_1 \cup S_2$, we have that $S_i \cap \pd N(f)$ is a single curve with slope $\frac{u}{1}$ for some $u$.  Hence, $S_i$ extends to a horizontal surface $S_i^*$ in the manifold $M^*$ which results from Dehn filling $C_{p,q} \setminus \eta(f)$ along the slope $\frac{u}{1}$.  Then $\chi(S_i^*) = \chi(S_i) + 1$ and $M^*$ is another cable space $C_{r,s}$.  Immediately, we have that $\chi(S_i^*) \leq -1$; thus $\chi(S) = \chi(S') = \chi(S_1) + \chi(S_2) \leq -4$. \\

Suppose further that the slope of each component of $S \cap \pd_+ C_{p,q}$ is $\frac{l}{1}$.  If $S_i^*$ is not planar, then $\chi(S_i^*) \leq -2$ and $\chi(S) \leq -6$ are desired, so suppose $S_i^*$ is planar.  Since the curves of $S \cap \pd_+ C_{p,q}$ intersect regular fibers of $C_{p,q}$ (which are also regular fibers of $C_{r,s}$) more than once, $S_i^*$ is not a horizontal planar surface of type 3(a) given by Lemma \ref{L1}.  Note that curves of slope $\frac{u}{1}$ and $\frac{1}{0}$ in $\pd_+ C_{p,q}$ will have respective slopes of $\frac{1}{0}$ and $\frac{v}{-1}$ for some integer $v$ in the Dehn filled manifold $C_{r,s}$.  The $2 \X 2$ matrix representing this change of basis is
\[ \begin{bmatrix} v & 1-uv \\ -1 & u \end{bmatrix}.\]
It follows that the slope of each component of $S^* \cap \pd_+ C_{r,s}$ is $\frac{1 - v(u-l)}{u-l}$, and by Lemma \ref{L1}, we have
\[ \frac{1-v(u-l)}{u-l} = \frac{1 + krs}{kr^2}.\]
Thus $-vkr^2 = -v(u-l) = krs$.  As $r$ and $s$ are relatively prime, we must have either $krs = 0$ or $krs = \pm 1$.  If $krs = 0$, then the slope of $S^* \cap \pd_- C_{r,s} = S \cap \pd_- C_{p,q}$ is meridional by Lemma 4 (since $C_{p,q}$ and $C_{r,s}$ have the same regular fibers).  On the other hand, if $krs = \pm 1$, then $kr^2 = \pm 1$, so $S^* \cap \pd_+ C_{r,s}$ has integral slope, implying $S^* \cap \pd_- C_{r,s}$ has integral slope by Lemma \ref{L2}.  This in turn shows that $S \cap \pd_- C_{p,q}$ has integral slope, completing the proof.

\end{proof}
\end{lemma}

For ease of notation but at the risk of confusion, we will let $C_0$ denote $E(K_0)$ despite the fact that $E(K_0)$ is not a cable space, and for $0 \leq l \leq m \leq n$, define
\[ C^m_l = C_l \cup C_{l+1} \cup \dots \cup C_m.\]
We prove one final lemma before the main theorem of this section.  

\begin{lemma}\label{L8}
Suppose that $\Sigma$ is a Heegaard surface for $C^m_l$ such that $\Sigma \cap C_l$ and $\Sigma \cap C_{l+1}$ are incompressible.  Then $\chi(\Sigma \cap (C_l \cup C_{l+1})) \leq -4$.
\begin{proof}
First, suppose that $S \cap C_l$ is a horizontal surface (this can only occur in the case $l = 0$), so each component of $\Sigma \cap C_l$ has negative Euler characteristic.  As $\Sigma$ is separating, there are at least 2 components; hence $\chi(\Sigma \cap C_l) \leq -2$.  In addition, $\Sigma \cap T_{l+1}$ has integral slope, so $\Sigma \cap C_{l+1}$ is horizontal and has at least 2 components, implying $\chi(\Sigma \cap C_{l+1}) \leq -2$ as well. \\

On the other hand, suppose that $\Sigma \cap C_l$ consists of vertical annuli.  By Lemma \ref{L1}, $\Sigma \cap C_{l+1}$ cannot be of type 3(a) since the regular fibers of $C_l$ and $C_{l+1}$ intersect more than once, and $\Sigma \cap C_{l+1}$ cannot be of type 3(b) since the slope of $\Sigma \cap T_{l+1}$ is integral.  It follows that each component $\Sigma_i$ of $\Sigma \cap C_{l+1}$ satisfies $\chi(\Sigma_i) \leq -2$; thus $\chi(
\Sigma \cap (C_l \cup C_{l+1})) = \chi(\Sigma \cap C_{l+1}) \leq -4$.

\end{proof}
\end{lemma}

\begin{theorem}\label{tunnel1}
For the iterated torus knot $K_n$, with cabling parameters satisfying $\Delta(\frac{p_i}{q_i}, p_{i-1}q_{i-1}) > 1$, we have $t(K_n) = n+1$.
\begin{proof}
First, we note that each $C_i$ has a minimal genus 2 Heegaard surface $\Sigma_i$, and together with the essential tori $T_i$, the collection $\{\Sigma_i\} \cup \{T_i\}$ yields a generalized Heegaard splitting of $E(K_n)$.  Amalgamating this splitting gives a Heegaard surface with genus $\sum g(\Sigma_i) - \sum g(T_i) = 2(n+1) - n = n+2$.  Thus, $t(K_n) \leq n+1$. \\

Fix $j$ such that $1 \leq j \leq n$, and suppose by way of induction that $g(C^{m'}_{l'}) = m'-l'+2$ whenever $m'-l' < j$.  The above argument covers the base case $m'-l' = 0$.  Let $\Sigma$ be a Heegaard surface for some $C^m_l$ such that $m-l = j$, and suppose first that $\Sigma$ is weakly reducible.  By Theorem \ref{thin}, untelescoping $\Sigma$ yields a generalized Heegaard splitting with at least one essential thin surface $S$, and by Lemma \ref{L5}, we have we have that $S = T_i$ for some $i$, where $l+1 \leq i \leq m$.  It follows that
\[ g(\Sigma) \geq g(C_l^{i-1}) + g(C_i^m) - g(T_i) = (i-1 - l + 2) + (m-i + 2) - 1 = m-l + 2,\]
completing the proof. \\

Hence, suppose that $\Sigma$ is strongly irreducible.  By Theorem \ref{graph}, there is an submanifold $M' = C_i$ or $T_i \X I$ of $C_l^m$ such that
\be
\item Each component of $\Sigma \setminus (\eta(\mathcal{T}) \cup \eta(M'))$ is incompressible in $M \setminus (\eta(\mathcal{T}) \cup \eta(M'))$.
\item If $M' = C_i$ for some $i$, then $\Sigma \cap M'$ is a pseudohorizontal surface in $M'$.
\item Otherwise, $M' = T_i \X I$, and $\Sigma \cap M'$ is either
{\be
\item the result of tubing an annulus parallel into $T_i \X \{0\}$ with boundary slope $s_0$ to an annulus parallel into $T_i \X \{1\}$ with boundary slope $s_1$ along a vertical arc in $T_i \X I$, where $\Delta(s_0,s_1) = 1$, or
\item the result of tubing two vertical annuli in $T_i \X I$ along an arc contained in $T_i \X \{\frac{1}{2}\}$.
\ee}
\ee

\noindent \textbf{Case 1}:  $M' = C_i$ for some $i$. \\

By the above, $\Sigma \cap C_i$ is pseudohorizontal and $\Sigma \cap C_j$ is incompressible for $j \neq i$.  Observe that $i \neq m$, and if $i = l$, then $l = 0$, as a pseudohorizontal surface intersects all boundary components of $C_i$.  Hence $\Sigma \cap C_m$ is a vertical annulus and $\chi(\Sigma \cap C_m) = 0$. \\

We will show that $\chi(\Sigma) \leq -2(m-l) - 2$.  If $\Sigma \cap C_j$ is a vertical annulus where $i < j < m$, then $\Sigma \cap \pd_+ C_{j+1}$ has integral slope.  By Lemma \ref{L2}, $\Sigma \cap \pd_- C_m \neq \emp$, contradicting that $\Sigma$ is a Heegaard surface for $C^m_l$.  If $i \neq l$, then $\Sigma \cap C_l$ is either a vertical annulus or a Seifert surface for $K_0$ (if $l=0$), so by Lemma \ref{L2}, $\Sigma \cap T_j$ has integral slope for $1\leq j \leq i$; thus $\Sigma \cap C_j$ is horizontal.  Thus, we may suppose that $\Sigma \cap C_j$ is horizontal whenever $j \neq i,l,m$.  As $\Sigma$ is separating, $\Sigma \cap C_j$ has at least two components, so $\chi(\Sigma \cap C_j) \leq -2$. \\

Suppose first that $i = 0$.  By Lemma \ref{L6}, we have
\begin{eqnarray*}
\chi(\Sigma) &=& \chi(\Sigma \cap C_0) + \chi(\Sigma \cap C_1) + \dots + \chi(\Sigma \cap C_{m-1}) + \chi(\Sigma \cap C_m) \\
&\leq& -4 - 2 - \dots - 2 - 0 \\
&\leq& -2(m-l) - 2.
\end{eqnarray*}

Next, suppose $i = l+1$.  By the remark above, $\Sigma \cap C_l$ is incompressible.  If each component of $\Sigma \cap C_l$ is horizontal, then $\chi(\Sigma \cap C_l) \leq -2$ and by Lemma \ref{L7}, $\chi(\Sigma \cap C_{l+1}) \leq -4$.  If $\Sigma \cap C_l$ is vertical annuli, then Lemma \ref{L7} provides that $\chi(\Sigma \cap C_{l+1}) \leq -6$ (or else $\Sigma \cap \pd_- C_m \neq \emp$ by Lemmas \ref{L2} and \ref{L3}, a contradiction).  In either case, we have
\begin{eqnarray*}
\chi(\Sigma) &=& \chi(\Sigma \cap (C_l \cup C_{l+1})) + \sum_{j=l+2}^{m-1} \chi(\Sigma \cap C_j) \leq -6 - 2(m-l-2) \\
&=& -2(m-l)-2.
\end{eqnarray*}

Finally, suppose that $i > l+1$.  Then $\Sigma \cap C_l$ and $\Sigma \cap C_{l+1}$ are incompressible, so by Lemma \ref{L8}, $\chi(\Sigma \cap (C_l \cup C_{l+1})) \leq -4$.  In addition, by Lemma \ref{L7}, $\chi(\Sigma \cap C_i) \leq -4$.  Otherwise, for $j \neq l,l+1,i,m$, components of $\Sigma \cap C_j$ are horizontal and $\Sigma \cap C_j$ contributes at most $-2$ to $\chi(\Sigma)$.  Thus
\begin{eqnarray*}
\chi(\Sigma) &\leq& \chi(\Sigma \cap (C_l \cup C_{l+1})) + \sum_{i=l+2}^{m-1} \chi(\Sigma \cap C_j) \\
&\leq& -4 - 2(m-l-2)-2 \\
&=& -2(m-l)-2.
\end{eqnarray*}

\noindent \textbf{Case 2}:  $M' = T_i \X I$ for some $i$. \\

Suppose that $T_i \X I \subset C_i$, so that $T_i = T_i \X \{0\}$ and let $T_i' = T_i \X \{1\}$.  In addition, set $C_i' = C_i \setminus \eta(T_i \X I)$ and for $j \neq i$, set $C_j' = C_j$, so that
\[ C^m_l = C_l' \cup_{T_{l+1}} \dots \cup_{T_{i-1}} C_{i-1}' \cup_{T_i} (T_i \X I) \cup_{T_i'} C_i' \cup_{T_{i+1}} \dots \cup_{T_m} C_m'.\]
By Theorem \ref{graph}, we have that $\chi(\Sigma \cap (T_i \X I)) = -2$, $\Delta((\Sigma \cap T_i),(\Sigma \cap T_i'))$ is either 0 or 1, and $\Sigma \cap C_j'$ is incompressible for all $j$.  By the argument above, each $\Sigma \cap C_j'$ must be a horizontal surface for $l \leq j < m$, with the possible exception of $\Sigma \cap C_l'$. \\

Suppose first that $i = l+1$.  If $\Sigma \cap C_l'$ is horizontal, then we have $\chi(\Sigma \cap C_j') \leq -2$ for all $j$, hence
\[ \chi(\Sigma) \leq \chi(\Sigma \cap (T_{l+1} \X I)) + \sum_{j=l}^{m-1} \chi(\Sigma \cap C_j') \leq -2 -2(m-l).\]
On the other hand, suppose that $\Sigma \cap C_l'$ is vertical, so that the slope of $\Sigma \cap T_{l+1}$ is $p_lq_l$.  If the slope of $\Sigma \cap T_{l+1}'$ is also $p_lq_l$, then by Lemma 2, $\Sigma \cap \pd_- C_m \neq \emp$, a contradiction.  It follows that the slope of $\Sigma \cap T_{l+1}'$ is $\frac{r}{s}$, where $r - p_lq_ls = \pm 1$. \\

Observe that $\Sigma \cap C_{l+1}'$ has at least two components, and if these are not planar, $\chi(\Sigma \cap C_{l+1}') \leq -4$, and again we have $\chi(\Sigma) \leq -2 - 2(m-l)$.  If $\Sigma \cap C_{l+1}'$ is planar, it must consist of two components each having one boundary component on $T_{l+1}'$.  By Lemma \ref{L1}, it follows that $\frac{r}{s} = \frac{1 + kp_{l+1}q_{l+1}}{kq_{l+1}^2}$.  Thus $1 + kp_{l+1}q_{l+1}- kp_{l}q_{l}q_{l+1}^2 = \pm 1$, and rearranging yields
\[ kq_{l+1}(p_{l+1} - p_lq_lq_{l+1}) = \pm 1 - 1.\]
By assumption $q_{l+1} \geq 2$ and $|p_{l+1} - p_lq_lq_{l+1}| > 1$, and so we must have $k=0$.  Therefore $\Sigma \cap C_i'$ is meridional, and by Lemma \ref{L3}, $\Sigma \cap \pd_- C^m \neq \emp$, a contradiction. \\

Finally, suppose $i > l+1$.  Then $\Sigma \cap (C_l \cup C_{l+1}) = \Sigma \cap (C_l' \cup C_{l+1}')$ and by Lemma \ref{L8}, $\chi(\Sigma \cap (C_l' \cup C_{l+1}')) \leq -4$.  Hence
\begin{eqnarray*}
\chi(\Sigma) &\leq& \chi(\Sigma \cap (C_l' \cup C_{l+1}')) + \chi(\Sigma \cap (T_i \X I)) + \sum_{j=l+2}^{m-1} \chi(\Sigma \cap C_j')\\
&\leq& -4 -2 -2(m-l-2) \\
&=& -2(m-l)-2,
\end{eqnarray*}
as desired.

\end{proof}
\end{theorem}

\section{Strongly irreducible bridge surfaces}\label{strong}

One valuable feature of strongly irreducible Heegaard surfaces for 3-manifolds is a ``no-nesting" property, demonstrated by Scharlemann in \cite{scharlemann}.  We adapt the proof of this important lemma to show a version of no-nesting for bridge surfaces below.

\begin{lemma}\label{L9}
Suppose $\Sigma$ is a strongly irreducible bridge surface in $(M,J)$, where $(M,J) = (V,\A) \cup_{\Sigma} (W,\n)$.
\be
\item If $c$ is an essential curve in $\Sigma_J$ such that $c$ bounds a disk $D \subset M(J)$, where a collar of $c$ in $D$ is disjoint from $\Sigma_J$, then $c$ bounds a compressing disk $D'$ in $(V,\A)$ or $(W,\n)$.
\item If $\gamma$ is an essential arc in $\Sigma_J$ such that $\gamma$ cobounds a disk $\Delta \subset M(J)$ with an essential arc $\mu \subset \pd N(J)$, where a collar of $\gamma$ in $\Delta$ is disjoint from $\Sigma_J$, then $\gamma$ cobounds a bridge disk in $(V,\A)$ or $(W,\n)$.
\ee
\begin{proof}
\be
\item Choose a disk $D$ with $\pd D = c$ and such that $|D \cap \Sigma_J|$ is minimal.  If $\text{int}(D) \cap \Sigma_J = \emp$, we are done.  If not, $D$ intersects $\Sigma_J$ in some number of simple closed curves.  Any curves which are inessential in $\Sigma_J$ may be removed by isotopy; hence, we may suppose that each curve of $D \cap \Sigma_J$ which is innermost curve in $D$ bounds a compressing disk for $\Sigma_J$. \\

Pick an innermost nested pair of curves $\delta$ and $\eps$, so that $\eps$ is innermost in $D$ and $\delta$ cobounds a component $P$ of $D \setminus \eta(\Sigma_J)$ with $\eps$ and possibly some other curves $\eps_1,\dots,\eps_n$, each of which is innermost in $D$ and thus bounds a compressing disk for $\Sigma_J$.  Suppose without loss of generality all of these disks lie in $(W,\n)$.  Compressing $\Sigma$ along this collection of disks yields a surface $\Sigma'$ and a bridge splitting $(M',J') = (V,\A) \cup_{\Sigma} (C,\n')$, where $\Sigma' = \pd_- C$, and $C$ is chosen so that $D \cap C$ contains a vertical annulus $A$ cobounded by $\delta$ and a curve $\delta' \subset D \cap \Sigma'$ (or in the case that $\delta = c$, the disk $D$ can be extended to $\Sigma'$ via such a vertical annulus $A$). \\

Clearly, $\delta'$ bounds a disk in $M'(J')$.  If $\delta'$ is essential in $\Sigma'_{J'}$, then by Theorem \ref{thin}, $\Sigma$ cannot be strongly irreducible, a contradiction.  Thus, $\delta'$ bounds a disk $D' \subset \Sigma'_{J'}$, and in $(M,J)$ we may perform surgery on $D$ by gluing $D'$ along $\delta'$.  However, this reduces the number of intersections of $\text{int}(D)$ with $\Sigma_J$, another contradiction.  We conclude that $\text{int}(D) \cap \Sigma_J = \emp$, completing the first part of the proof. \\

\item The second statement follows easily by viewing $\Delta$ as a disk in $M$ with a boundary arc contained in $J$ and by observing that $N(\Delta) \subset M$ is a 3-ball, where $c = \pd (N(\gamma) \cap \Sigma)$ bounds a compressing disk $D \subset \pd N(\Delta)$ satisfying (1).  Thus, $c$ bounds a disk $D'$ in $(V,\A)$ or $(W,\n)$.  It must then be true that the arc $\mu \subset \pd N(J)$ intersects $\Sigma_J$ only in its endpoints; hence $D'$ is the frontier in $V$ or $W$ of a regular neighborhood of a bridge disk $\Delta'$, where $\gamma \subset \pd \Delta'$.

\ee

\end{proof}
\end{lemma}

Next, we adapt a lemma from the theory of Heegaard surfaces in the context of bridge surfaces \cite{bachschsedg}, \cite{li}.  The lemma asserts that strongly irreducible surfaces behave much like incompressible surfaces; namely, cutting a strongly irreducible surface $A$ along a collection of incompressible surfaces splits $A$ into a number of pieces, all of which are incompressible with the exception of at most one strongly irreducible component.  We make this statement rigorous in the next lemma, whose proof is modeled on the proof of Lemma 3.7 in \cite{li}.  We need several definitions before proceeding. \\

First, we weaken the definition of an essential surface.  Let $M$ be a 3-manifold with boundary, $P \subset \pd M$ a subsurface, and $A \subset M$ a properly embedded surface.  A $P$-$\pd$-\emph{compressing} disk for $A$ is a $\pd$-compressing disk $\Delta$ for $A$ such that $\Delta \cap \pd M \subset P$.  We say that $A$ is $P$\emph{-essential} if $A$ is incompressible and there does not exist a $P$-$\pd$-compressing for $A$ in $M$.  On the other hand, if $A$ is separating and admits compressing or $P$-$\pd$-compressing disks on either side but admits no pair of disjoint disks on opposite sides, we say that $A$ is $P$\emph{-strongly irreducible}.  As an example, a strongly irreducible bridge surface for a knot $K$ in $M$ is $\pd N(K)$-strongly irreducible in $M(K)$.  Finally, we say that two surfaces $A$ and $B$ are \emph{almost tangent} if $A$ is transverse to $B$ except for a single saddle tangency. \\

The lemma below is stated in the greater generality than is needed here; however, we include it in its entirety for anticipated use in future work.

\begin{lemma}\label{L10}
Let $M$ be a compact 3-manifold and $J$ a properly embedded 1-manifold, with $Q = \pd N(J)$ in $M(J)$.  Suppose $\Sigma$ is a strongly irreducible bridge splitting surface for $(M,J)$, and let $S \subset M(J)$ be a collection of properly embedded essential surfaces such that for each component $c$ of $\pd S$, either $c \subset Q$ or $c \subset \pd M$.  Then one of the following must hold:

\be
\item[(1)] After isotopy, $\Sigma_J$ is transverse to $S$ and each component of $\Sigma_J \setminus \eta(S)$ is $Q$-essential in $M(J) \setminus \eta(S)$.

\item[(2)] After isotopy, $\Sigma_J$ is transverse to $S$, one component of $\Sigma_J \setminus \eta(S)$ is $Q$-strongly irreducible and all other components are $Q$-essential in $M(J) \setminus \eta(S)$,

\item[(3)] After isotopy, $\Sigma_J$ is almost transverse to $S$, and each component of $\Sigma_J \setminus \eta(S)$ is $Q$-essential in $M(J) \setminus \eta(S)$.

\ee

\begin{proof}
Let $(M,J) = (V,\A) \cup_{\Sigma} (W,\n)$ be a bridge splitting with strongly irreducible bridge surface $\Sigma$, and let $G_V$ and $G_W$ denote cores of $V$ and $W$, respectively.  If $V = \Sigma \X I$, define $G_V$ to be a small arc with one endpoint on $\pd_- V$ and one endpoint in $\text{int}(V)$ (similarly with $W$ and $G_W$).  After isotopy, we may assume that $G_V$ and $G_W$ miss each arc in $\A$ and $\n$.  Let $b_{\A}$ be the number of $\pd$-parallel arcs in $(V,\A)$, and define $\Gamma_V$ to be the union of $G_V$, $(\pd_- V)_J$, $b_{\A}$ meridional curves contained in $\pd N(\A) \subset V(\A)$, and $b_{\A}$ unknotted arcs connecting the meridional curves to $G_V$, so that $V(\A) \setminus \Gamma_V \cong \Sigma_J \X (0,1]$.  We may define $\Gamma_W \subset W(\n)$ similarly. \\

Then $M(J) \setminus (\Gamma_V \cup \Gamma_W) \cong \Sigma_J \X (0,1)$.  This induces a sweepout $f:M(J) \rightarrow I$ such that $f^{-1}(0) = \Gamma_V$, $f^{-1}(1) = \Gamma_W$, and for every $t \in (0,1)$, $\Sigma_t = f^{-1}(t)$ is isotopic to $\Sigma_J$.  We let $(M,J) = (V_t,\A_t) \cup_{\Sigma_t} (W_t,\n_t)$ denote the bridge splitting induced by $\Sigma_t$.  Assume that $S$ is transverse to $\Gamma_V \cup \Gamma_W$.  Choosing $\eps$ small enough, we then have that both $S \cap ((V_{\eps},\A_{\eps}) \cup (W_{1-\eps},\n_{1-\eps}))$ is a collection of compressing and bridge disks. \\

After a small perturbation, the function $h= f|_S$ is Morse on $(-1,1)$, where all critical points of $h$ occur at different levels, each meridional component of $\pd S$ is contained in a unique level away from critical values of $h$ (these levels will also be considered critical values, viewed as truncated extrema), and for any other component $\gamma$ of $\pd S$, $h|_{\gamma}$ is monotone on $(-1,1)$, so that for any two components $c_1 \subset \pd \Sigma_t$ and $c_2 \subset \pd S$, we have $|c_1 \cap c_2|$ is minimal us to isotopy. \\

We assign to each $t \in (-1,1)$ some subset of the labels $\{\nu,\omega\}$ by the following method:  If $\Sigma_t \setminus \eta(S)$ contains a curve which bounds a compressing disk in $(V_t,\A_t)$ or an arc which cobounds a bridge disk in $(V_t,\A_t)$, we assign $t$ the label $\nu$.  The label $\omega$ is defined analogously with $(W_t,\n_t)$.  Note that $t$ may have either, both, or no labels.  In addition, we observe that $\eps$ has the label $\nu$ and $1-\eps$ has the label $\omega$, which follows from the fact that $S$ is transverse to $\Gamma_V$ and $\Gamma_W$, and both cores are nonempty. \\

Next, we claim that if some regular value $t$ of $h$ has no label, then (1) holds.  Suppose $t$ has no label.  We assert that each curve and arc of $S \cap \Sigma_t$ is essential in both $S$ and $\Sigma_t$ or inessential in $S$ and $\Sigma_t$.  By the incompressibility and $\pd$-incompressibility of $S$, no curve or arc is essential in $S$ but inessential in $\Sigma_t$.  Suppose that $S \cap \Sigma_t$ contains a curve or arc $c$ which is essential in $\Sigma_t$ but inessential in $S$.  By the no-nesting Lemma above and the strong irreducibility of $\Sigma_t$, this implies that $c$ bounds a compressing disk or cobounds a bridge disk in either $(V_t,\A_t)$ or $(W_t,\n_t)$, contradicting the assumption that $\Sigma_t$ has no label.  Thus, curves and arcs in $S \cap \Sigma_t$ are either essential or inessential in both. \\

If $S \cap \Sigma_t$ contains an arc $\gamma$ which is inessential in both $\Sigma_t$ and $S$, then both endpoints of $\gamma$ are contained in single component $c_1$ of $\pd \Sigma_t$ and $c_2$ of $\pd S$.  However, this implies that $\pd \gamma \subset c_1 \cap c_2$ contains two points of opposite algebraic intersection number, contradicting that $|c_1 \cap c_2|$ is minimal.  Thus, arcs of $S \cap \Sigma_t$ are essential in both surfaces. \\

Let $c$ be a curve in $S \cap \Sigma_t$ which is innermost among inessential curves in $S$.  Then $c$ bounds a disk $\Delta \subset S$ which misses $\Sigma_t$, and performing surgery on $\Sigma_t$ along $c$ gives a surface isotopic to $\Sigma_t$ with fewer intersections with $S$.  Finitely many iterations of this operation yields a bridge surface $\Sigma_t'$ isotopic to $\Sigma_t$ and such that each curve and arc in $\Sigma_t' \cap S$ is essential in both surfaces.  Moreover, $\Sigma_t' \cap S \subset \Sigma_t \cap S$. \\

Suppose now that some component of $\Sigma_t' \setminus \eta(S)$ is $Q$-compressible.  Then there is a curve or arc $c \subset \Sigma_t' \setminus \eta(S)$ such that $c$ bounds a disk or cobounds a bridge disk in $(V_t,\A_t)$ or $(W_t,\n_t)$.  Since each curve of $\Sigma_t' \cap S$ is essential in $\Sigma_t'$, it follows that $c$ is essential in $\Sigma_t'$.  We may isotope $c$ so that $c \subset \Sigma_t$, but this means that $t$ is labeled $\nu$ or $\omega$, a contradiction.  We conclude that each component of $\Sigma_t' \setminus \eta(S)$ is $Q$-essential in its respective submanifold. \\

Our second major claim is that if a regular value $t$ is labeled $\nu$ and $\omega$, then (2) holds.  As above, we assert that $S \cap \Sigma_t$ contains curves and arcs that are either essential in both $S$ and $\Sigma_t$ or inessential in both $S$ and $\Sigma_t$.  The incompressibility and $\pd$-incompressibility of $S$ rules out curves and arcs essential in $S$ but inessential in $\Sigma_t$.  Suppose $c$ is a curve or arc in $S \cap \Sigma_t$ which is essential in $\Sigma_t$ but inessential in $S$.  Then $c$ bounds or cobounds a disk in $M(J)$, and by the no-nesting Lemma above, $c$ bounds or cobounds a disk in $(V_t,\A_t)$ or $(W_t,\n_t)$.  However, since $t$ is labeled $\nu$ and $\omega$, there exist curves or arcs $c_V$ and $c_W$ disjoint from $S$ (and thus from $c$) which bound or cobound disks in $(V_t,\A_t)$ and $(W_t,\n_t)$, respectively.  This contradicts the strong irreducibility of $\Sigma_t$.  Hence, curves and arcs in $S \cap \Sigma_t$ are either essential or inessential in both surfaces. \\

As above, $S \cap \Sigma_t $ cannot contain arcs inessential in both surfaces, and we may construct a surface $\Sigma_t'$ such that $S \cap \Sigma_t'$ contains only curves and arcs which are essential in both $S$ and $\Sigma_t'$.  Let $\Sigma'$ be the component of $\Sigma_t' \setminus \eta(S)$ which contains $c_V$.  By the strong irreducibility of $\Sigma_t'$, $c_V \cap c_W \neq \emp$ and thus $c_W \subset \Sigma'$ and, more generally, $\Sigma'$ must be $Q$-strongly irreducible.  Let $\Sigma''$ be any other component of $\Sigma_t' \setminus \eta(S)$.  If an essential curve or arc $c$ bounds or cobounds a disk for $(V_t,\A_t)$ or $(W_t,\n_t)$, this gives rise to a compressing or bridge disk for $\Sigma_t$ disjoint from $c_V$ and $c_W$, another contradiction.  We conclude that $\Sigma''$ is $Q$-essential in $M(J) \setminus \eta(S)$. \\

In the final remaining case, suppose that $h$ has a critical value $t \in (\eps,1-\eps)$ such that $t - \delta$ is labeled $\nu$ and $t + \delta$ is labeled $\omega$, with curves or arcs $c_V$ and $c_W$ contained in $\Sigma_{t\pm\delta} \setminus \eta(S)$ bounding or cobounding disks in $(V_{t-\delta},\A_{t-\delta})$ and $(W_{t+\delta},\n_{t+\delta})$, respectively.  If $c$ corresponds to minimum, maximum, or a level component of $\pd S$, then there exists an isotopy pushing $c_V$ to $c'_V \in \Sigma_{t + \delta}$ or $c_W$ to $c_W'$ in $\Sigma_{t-\delta}$ bounding disks in $(V_{t+\delta},\A_{t+\delta})$ or $(W_{t-\delta},\n_{t-\delta})$, respectively, contradicting either the assumption that $t + \delta$ is labeled $\omega$ or the assumption that $t-\delta$ is labeled $\nu$.  We conclude that $t$ corresponds to a saddle. \\

We may regard a small closed regular neighborhood $N(S)$ of in $M(J)$ as $S \X I$.  Then $\Sigma_t \cap (S \X I)$ contains components of the form $\gamma \X I$, where $\gamma$ is curve or arc in $\Sigma_t \cap S$, in addition to one exceptional component $E$ containing the saddle point.  The surface $E$ must be a pair of pants, an annulus or a disk, depending on whether $S \cap \Sigma_t$ is a figure-8 curve, a curve wedged with an arc, or the wedge of two arcs.  Figure \ref{sads} depicts the possible configurations corresponding to the saddle point along with potential singular sets $S \cap \Sigma_t$.  Note while the figure depicts subsets of $S$ on which $h$ is a Morse function, each of these pieces is isotopic to the corresponding exceptional component $E \subset \Sigma_t$. \\

\begin{figure}[h!]
  \centering
    \includegraphics[width=0.8\textwidth]{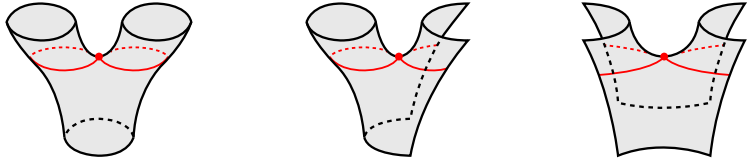}
    \caption{Possible components of $S$ corresponding to a saddle}
    \label{sads}
\end{figure}

We assert that if a curve or arc $\gamma \subset \Sigma_t \setminus \eta(S)$ is essential in $\Sigma_t$, then $\gamma$ does not bound or cobound a disk in $M$.  Otherwise, by the no-nesting Lemma $\gamma$ bounds or cobounds a disk in $(V_t,\A_t)$ or $(W_t,\n_t)$, implying that $\Sigma_{t - \delta}$ and $\Sigma_{t + \delta}$ have the a common label, since for small $\delta$, $\Sigma_{t\pm \delta} \setminus \eta(S)$ is parallel to $\Sigma_t \setminus \eta(S)$. \\

Next, we claim that curves and arcs in $\Sigma_t \cap \pd N(S)$ are either essential or inessential in both $\pd N(S)$ and $\Sigma_t$.  Since $\pd N(S)$ is essential, no curve or arc $\gamma \in \pd N(S) \cap \Sigma_t$ is essential in $\pd N(S)$ and inessential in $\Sigma_t$.  Additionally, by the assertion above, no such $\gamma$ is essential in $\Sigma_t$ and inessential in $\pd N(S)$.  As above, there cannot be arcs in $\pd N(S) \cap \Sigma_t$ which are inessential in both surfaces.  If there exists a curve of intersection which is inessential in both surfaces, cutting and pasting along a curve which is innermost in $\pd N(S)$ yields a surface isotopic to $\Sigma_t$ with fewer intersections with $\pd N(S)$, and finitely many repetitions produces a surface $\Sigma_t'$ isotopic to $\Sigma_t$ and such that $\pd N(S) \cap \Sigma'_t$ contains only curves and arcs which are essential in both surfaces. \\

Here we show that every component of $\Sigma'_t \setminus \eta(S)$ is $Q$-essential.  If $\gamma$ is a curve or arc which is essential in $\Sigma'_t \setminus \eta(S)$, then $\gamma$ must be essential in $\Sigma'_t$ as all components of $\Sigma'_t \cap \pd N(S)$ are essential curves or arcs.  Thus, $\gamma$ cannot bound or cobound a disk in $M$.  Otherwise,  after isotopy we may assume $\gamma \subset \Sigma_t$ and by the no nesting Lemma, $\gamma$ bounds or cobounds a disk in $(V_t,\A_t)$ or $(W_t,\n_t)$, contradicting our assertion above. \\

To finish the proof, we must show that $\Sigma_t'$ is either tangent or almost tangent to $S$.  In process of capping off disks to get $\Sigma'_t$ from $\Sigma_t$, we may have altered the exceptional component $E$ if some curve in $\pd E$ was inessential in $\pd N(S)$.  Let $E' = E \cap \Sigma_t'$ (that is, $E'$ is the essential subsurface of $\Sigma_t'$ resulting from capping off inessential curves and arcs of $E$).  There are several cases to consider:  If $E' = \emp$, then (1) holds.  If $E' = E$ then (3) holds. \\

If $E' \neq \emp$ and $E' \neq E$, then the above process must have capped off exactly one curve component of $\pd E$, and $\pd E'$ has two components, $c_0$ and $c_1$.  If $c_0$ and $c_1$ are contained in different components of $\pd N(S)$, then $E'$ can be expressed as $c_0 \X I$ in $N(S)$, and (1) holds.  Otherwise, $c_0$ and $c_1$ are in the same component, call it $S_0$, of $\pd N(S)$.  We may replace $S$ with $S_0$ (since the two surfaces are isotopic), noting that $\Sigma_t'$ is transverse to $S_0$ and each component of $\Sigma_t' \setminus \eta(S_0)$ is $Q$-essential in its respective submanifold; that is, (1) holds, as desired.

%Let $S_0$ and $S_1$ denote components of $\pd N(S) = S \X \pd I$.  If $\pd E'$ contaiIf $\pd E'$ is contained entirely in, say, $S_0$, we may replace $S$ with $S_0$ (since the two surfaces are isotopic), noting that $\Sigma_t'$ is transverse to $S_0$ and each component of $\Sigma_t' \setminus \eta(S_0)$ is $J$-incompressible in its respective submanifold; that is, (1) holds.  

%To finish the proof, suppose that $E' \neq E$, $E' \neq \emp$, and $E'$ has boundary components on both $S_0$ and $S_1$.  This assumption yields two possibilities:  First, suppose that $E'$ can be expressed as $\gamma \X I \subset S \X I$ for a curve or arc $\gamma$.  Then after perturbation, $\Sigma_t'$ is transverse to $S$ and (1) holds.  If $E'$ cannot be expressed as $\gamma \X I$, we construct a new surface $S'$ by attaching $S_0 \cap V_t'$, $S_1 \cap W_t'$, and $(S \X I) \cap \Sigma_t'$ along their shared boundary components (the roles of $S_0$ or $S_1$ may need to be switched).  Then $S'$ is isotopic to $S$, and after a slight perturbation, $S'$ is transverse to $\Sigma'$ except for a half saddle point contained in $E'$.  In this case, (3) holds, completing the proof.

\end{proof}
\end{lemma}

\section{Analyzing the bridge spectra of $K_n$}\label{spec}

In Section \ref{tunnel}, we utilized Theorem \ref{graph} to understand minimal genus Heegaard surfaces for $E(K_n)$.  In a similar vein, in this section we will employ Lemma \ref{L10} to characterize minimal bridge surfaces for $K_n$, although the analysis here is significantly less complicated than that of Section \ref{tunnel}.  We begin with a lemma concerning embeddings of cables on Heegaard surfaces.

\begin{figure}[h!]
  \centering
    \includegraphics[width=1.0\textwidth]{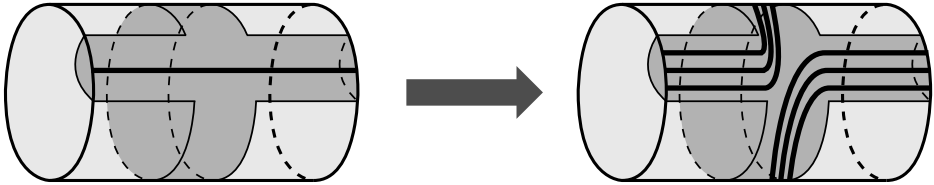}
    \caption{Replacing $K \subset \Sigma$ with its cable}
        \label{cabheg}
\end{figure}

\begin{lemma}\label{L11}
Let $K$ be a knot in a 3-manifold $M$ and $K_{p,q}$ a $(p,q)$-cable of $K$.  If $\Sigma \subset M$ is a Heegaard surface such that $K \subset \Sigma$ and $D$ is a compressing disk for $\Sigma$ such that $|D \cap K| = 1$, then there exists an embedding of $K_{p,q}$ in $M$ such that $K_{p,q} \subset \Sigma$.
\begin{proof}
Suppose $M = V \cup_{\Sigma} W$, with $D \subset V$, and let $D \X I$ be a collar neighborhood of $D$ in $V$, noting that $K \cap (\pd D \X I)$ is a single essential arc.  In addition, let $X = N(K) \cap V$.  Then both $X$ and $X' = X \cup (D \X I)$ are solid tori whose cores are isotopic to $K$.  Now $X' \cap S = \pd X' \cap S$ is a once-punctured torus.  Since every torus knot can be embedded on a once punctured torus, every $(p,q)$-cable of $K$ can be embedded on $X' \cap S$.  See Figure \ref{cabheg}.
\end{proof}
\end{lemma}

As a consequence of Lemma \ref{L11}, we may demonstrate that $b_g(K_n) = 0$ whenever $g > n$:  By Theorem \ref{tunnel}, $t(K_{n-1}) = n$, which implies that $E(K_{n-1})$ has a genus $n+1$ Heegaard splitting.  In this case, the $\infty$-sloped Dehn filling of $K_{n-1}$ yields a genus $n+1$ Heegaard surface $\Sigma$ of $S^3$ in which $K_{n-1}$ is contained in a core of one of the handlebodies cut out by $\Sigma$.  It follows that $K_{n-1}$ is isotopic into $\Sigma$, and there is a compressing disk $D$ such that $|D \cap K_{n-1}| = 1$.  By Lemma \ref{L11}, there is an embedding of $K_n$ such that $K_n \subset \Sigma$.  Hence $b_{n+1}(K_n) = 0$, and using stabilization, $b_g(K_n) = 0$ for all $g > n$. \\

Conversely, if $K_n$ is isotopic into a genus $g$ Heegaard surface $\Sigma'$, then a perturbation followed by a meridional stabilization of $\Sigma'$ yields a genus $g+1$ surface $\Sigma''$ such that there is a compressing disk $D'$ for $\Sigma''$ with $|D' \cap K_n| = 1$.  In this case, $K_n$ is isotopic into the core of a handlebody cut out by $\Sigma''$, implying that $g(E(K_n)) \leq g+1$.  This implies that $t(K_n) \geq g$; thus by Theorem \ref{tunnel}, $n+1$ must be the smallest $g$ for which $b_g(K_n) = 0$, and $\Sigma$ is not cancelable for $g \leq n$. \\

For any knot $K \subset S^3$, the smallest $g$ for which $b_g(K) = 0$ is called the \emph{h-genus} $h(K)$ of $K$, following Morimoto \cite{morim}.  In this work, Morimoto partitions the set of all knots into sets $A_n$, $B_n$, and $C_n$ related to tunnel number, h-genus, and another invariant known as 1-bridge genus.  He conjectures that each of these sets is nonempty, and the above argument verifies that for all $K_n$, we have $h(K_n) = t(K_n) = n+1$ and $K_n \subset A_{n+1}$. \\

Before we arrive at the proof of the main theorem, recall that $E(K_n)$ decomposes as
 \[ E(K_n) = C_0 \cup_{T_1} \cup \dots \cup_{T_n} C_n,\]
 and that $\mathcal{T}$ denotes the collection $\{T_1,\dots,T_n\}$.  Note also that $C_n$ contains an essential vertical annulus $A$ such that $\pd A \subset \pd E(K_n)$, where $C_n \setminus \eta(A)$ is the union of $T_n \X I$ and a solid torus $X$.  Thus, if we set $\mathcal{T}' = \mathcal{T} \cup \{A\}$, then $E(K_n)$ cut along $\mathcal{T'}$ is the union of $C_0$, the cable spaces $C_1,\dots,C_{n-1}$, the product region $T_n \X I$, and the solid torus $X$. \\
 
Define $V_0 = S^3$, and for every $l$ with $1 \leq l \leq n$, let $V_l$ denotes the solid torus in $S^3$ bounded by $T_l$, so that $V_l \setminus \eta(K_n) = C^n_l$ (using the terminology of Section \ref{tunnel}).

\begin{lemma}\label{L12}
Suppose that $\Sigma$ is a strongly irreducible $(g,b)$-bridge surface for $(V_l,K_n)$ with $0 \leq l \leq n-1$ and $b \geq 1$.  Then one of the following holds:
{\be
\item $\Sigma$ is a $(g,b')$-bridge surface for $(V_l,K_{n-1})$, where $b \geq q_n \cdot b'$ and $b' \geq 1$,
\item $g > n-l$, or
\item $l = n-1$, $g = 1$, and $b \geq |p_n - p_{n-1}q_{n-1}q_n|$.
\ee}
\begin{proof}
Let $Q = \pd N(K_n)$.  By Lemma \ref{L10}, there is an isotopy of $\Sigma$ after at most one of $\Sigma \cap C_l,\dots \Sigma \cap C_{n-1}, \Sigma \cap (T_n \X I), \Sigma \cap V$ is $Q$-strongly irreducible, while the remaining surfaces are $Q$-essential in their respective submanifolds.  Suppose first that $\Sigma \cap X$ is $Q$-essential.  Since $\Sigma \cap K_n \neq \emp$, we have that $\Sigma \cap X$ cannot be a $\pd$-parallel annulus or disk, as these surfaces are either disjoint from $K_n$ or are $Q$-$\pd$-compressible.  Thus each component of $\Sigma \cap X$ is a meridian disk, which intersects $K_n$ at least $q_n$ times.  This implies that $\Sigma$ is a $(g,b')$-bridge surface for a core $K_{n-1}$ of $X$ intersecting $\Sigma$ transversely, implying $b \geq q_n \cdot b'$ and $b' \geq 1$. \\

On the other hand, suppose that $\Sigma \cap X$ is $Q$-strongly irreducible.  Then $\Sigma \cap (T_n \X I)$ is incompressible, and thus must be a collection of vertical annuli.  This implies that $\Sigma \cap T_n$ is essential simple closed curves and parametrizing $\pd X$ as $T_n$, slopes of $\Sigma \cap \pd X$ and $\Sigma \cap T_n$ are equal.  In addition, each $\Sigma \cap C_j$ is essential.  It follows that $\Sigma \cap C_l$ must be a collection of vertical annuli or Seifert surfaces for $K_0$ (if $l=0$); thus by Lemma \ref{L2}, $\Sigma \cap T_j$ has integral slope for all $j$.  Thus, $\Sigma \cap X$ has no disk components and $\chi(\Sigma \cap X) \leq 0$.  If $l \leq n-2$, then by Lemma \ref{L8} $\chi(\Sigma \cap (C_l \cup C_{l+1})) \leq -4$ and $\chi(\Sigma \cap C_j) \leq -2$ whenever $l+1 < j < n$.  Therefore
\[ \chi(\Sigma) = \chi(\Sigma \cap (C_l \cup C_{l+1})) + \sum_{j=l+2}^{n-1} \chi(\Sigma \cap C_j) \leq -4 - 2(n-l-2),\]
$g(\Sigma) \geq n-l + 1$, and (2) holds. \\

Thus, if $l \leq n-2$, then (1) or (2) holds.  Now suppose $l = n-1$, so if (2) does not hold, then $\Sigma$ is a $(1,b)$-surface.  As above, if $\Sigma \cap X$ is meridian disks, (1) holds.  Otherwise, $\Sigma \cap C_{n-1}$ is a vertical annulus, as is $\Sigma \cap (T_n \X I)$.  Viewing $K_n$ as a regular fiber of $C_n$, we see that each component of $\Sigma \cap T_n$ must contribute $|p_n - p_{n-1}q_{n-1}q_n|$ intersections with $K_n$, and since there are at least two such components, we have $b \geq |p_n - p_{n-1}q_{n-1}q_n|$.

\end{proof}
\end{lemma}

Note that by \cite{schub} and \cite{schult}, the bridge spectrum of the $(p_0,q_0)$-torus knot $K_0$ is $\mathbf{b}(K_0) = (\min\{p_0,q_0\},0)$.  We describe the spectrum of $K_n$ inductively:
  
 \begin{theorem}\label{spectra}
 Suppose that $K_n$ is an iterated torus knot, whose cabling parameters satisfy $|p_i - p_{i-1}q_{i-1}q_i| > 1$.  Then
 \[b_g(K_n)=	\begin{cases}
			q_n \cdot b_g(K_{n-1}) &\text{if $g < n$;}\\
			\min\{|p_n - p_{n-1}q_{n-1}q_n|,q_n\} &\text{if $g = n$;}\\
			0 &\text{otherwise.}
		\end{cases}
\]
In other words,
\[ \mathbf{b}(K_n) = q_n \cdot \mathbf{b}(K_{n-1}) + \min\{|p_n - p_{n-1}q_{n-1}q_n|,q_n\} \cdot \mathbf{e}_n.\]

\begin{proof}
By Theorem \ref{satbridge}, $b_0(K_n) = q_n \cdot b_0(K_{n-1})$, and by the above arguments $b_g(K_n) = 0$ if and only if $g > n$.  First, we exhibit bridge surfaces satisfying the above equalities.  Any $(g,b)$-bridge surface $\Sigma$ for $K_{n-1}$ can also be seen as a $(g,q_n \cdot b)$-surface for $K_n$ by replacing each trivial arc of $K_{n-1}$ with $q_n$ trivial arcs of its cable $K_n$.  If follows that $b_g(K_n) \leq q_n \cdot b_g(K_{n-1})$ for $g < n$, and since $b_n(K_{n-1}) = 0$, perturbing this surface yields an $(n,1)$-surface for $K_{n-1}$, implying $b_n(K_n) \leq q_n$.  \\

Let $p_n^* = |p_n - p_{n-1}q_{n-1}q_n|$.  Following the proof of Lemma \ref{L12}, we can see that $(V_{n-1},K_n)$ has a $(1,p_n^*)$-surface $\Sigma'$ constructed by taking the union of a vertical annulus in $C_{n-1}$, a vertical annulus in $T_n \X I$, and a $\pd$-parallel annulus in $X$.  If $n = 1$, this is a $(1,p_n^*)$-surface for $(S^3,K_1)$; otherwise, by Theorem \ref{tunnel1}, $E(K_{n-2})$ has a genus $n$ Heegaard surface $\Sigma''$, and thus $\{\Sigma',T_{n-1},\Sigma''\}$ is a generalized bridge splitting for $(S^3,K_n)$ whose amalgamation $\Sigma^*$ is an $(n,p_n^*)$-bridge surface.  It follows that $b_n(K_n) \leq p_n^*$. \\

Now, let $(S^3,K_n) = (V,\A) \cup_{\Sigma} (W,\n)$ be a $(g,b)$-bridge splitting with $g \leq n$, so that $b \geq 1$ by above arguments.  If $\Sigma$ is strongly irreducible, then by Lemma \ref{L12}, either $b \geq q_n \cdot b_g(K_{n-1})$, or $g = n = 1$ and $b \geq \min\{p_n^*,q_n\}$.  Otherwise, $\Sigma$ is weakly reducible and may be untelescoped to yield a generalized bridge splitting $\{\Sigma_0,S_1,\dots,S_d,\Sigma_d\}$, where each $\Sigma_i$ is strongly irreducible and each $S_i$ is essential.  By Theorem 6.6 of \cite{zupan}, $E(K_n)$ contains no essential meridional surfaces, so $S_i \cap K_n = \emp$ for all $i$ and thus by Lemma \ref{L5}, each $S_i$ must be isotopic to some $T_j$. \\

As $K_n$ lies on the same side of $T_j$ for all $j$, either $\Sigma_0 \cap K_n \neq \emp$ or $\Sigma_d \cap K_n \neq \emp$.  Assume $\Sigma_d \cap K_n \neq \emp$ and consider let $T_l = S_d$.  Then $\Sigma_d$ is a strongly irreducible bridge surface for some $V_l$, $\Sigma_d$ has the same bridge number as $\Sigma$, and $\{\Sigma_0,S_1,\dots,S_{d-1},\Sigma_{d-1}\}$ is a generalized Heegaard splitting for $E(K_{l-2})$.  By Theorem \ref{tunnel1}, we have
\[ \sum_{i = 0}^{d-1} g(\Sigma_i) - \sum_{i=1}^{d-1} g(S_i) \geq g(E(K_{l-2})) = l+1.\]
Consider the three possibilities afforded to $\Sigma_d$ by Lemma \ref{L12}:  In the first case, $b \geq \max\{q_n \cdot b_g(K_{n-1}),q_n\}$.  If $g(\Sigma_d) > n-l$, then we have
\[ g(\Sigma) = g(\Sigma_d) - g(T_l) + \sum_{i = 0}^{d-1} g(\Sigma_i) - \sum_{i=1}^{d-1} g(S_i) > (n -l) - 1 + (l+1) = n,\]
a contradiction. \\

In the final case, $l = n-1$, $g(\Sigma_d) = 1$, and $b \geq p_n^*$.  This implies
\[ g(\Sigma) = g(\Sigma_d) - g(T_{n-1}) + \sum_{i = 0}^{d-1} g(\Sigma_i) - \sum_{i=1}^{d-1} g(S_i) \geq 1 - 1 + n = n;\]
hence $g(\Sigma) = n$ and $b_n(K_n) \geq p_n^*$, completing the proof of the theorem.

\end{proof}
\end{theorem}
 
\section{An example}\label{cableex}

\begin{figure}[h!]
  \centering
    \includegraphics[width=0.8\textwidth]{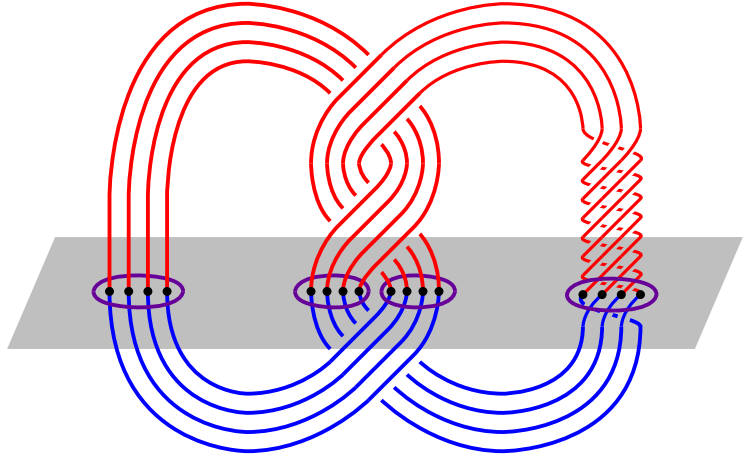}
    \caption{A $(0,8)$-surface $\Sigma_0$ for $K_1$, with purple curves depicting $T_n \cap \Sigma_0$}
        \label{bridge0}
\end{figure}

Here we apply Theorem \ref{spectra} to produce the bridge spectrum of $K_1 = ((3,2),(21,4))$, with illustrations.  Note that $K_1$ is a $(21,4)$-cable of the trefoil $K_0 = ((3,2))$, so by Theorem \ref{satbridge}, we have $b_0(K_1) = 8$.  An illustration of a minimal $(0,8)$-surface $\Sigma_0$ appears in Figure \ref{bridge0}, where $X \cap \Sigma_0$ is a collection of meridian disks, each intersecting $K_1$ four times. \\

\begin{figure}[h!]
  \centering
    \includegraphics[width=0.65\textwidth]{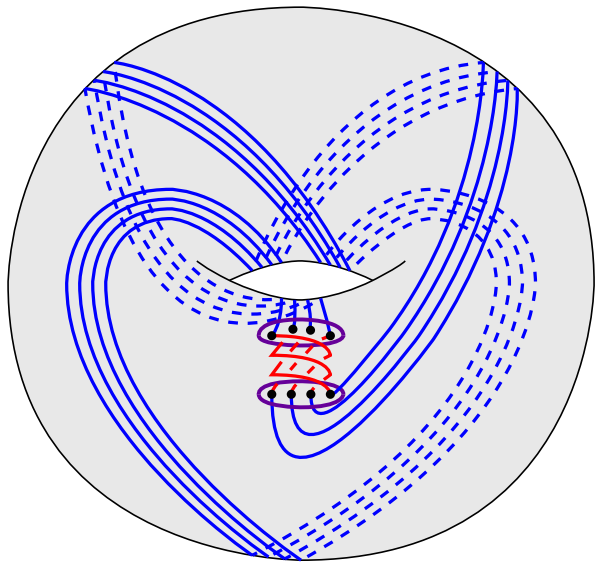}
    \caption{A $(1,4)$-surface $\Sigma_1$ for $K_1$, with purple curves depicting $T_n \cap \Sigma_1$}
        \label{bridge1}
\end{figure}

Turning to the genus one case, taking an obvious cabling of $(1,1)$-surface for $K_0$ yields a $(1,4)$-surface $\Sigma_1$ for $K_1$, where $X \cap \Sigma_1$ is a collection of meridian disks, each of which hits $K_1$ four times.  See Figure \ref{bridge1}. \\

\begin{figure}[h!]
  \centering
    \includegraphics[width=0.65\textwidth]{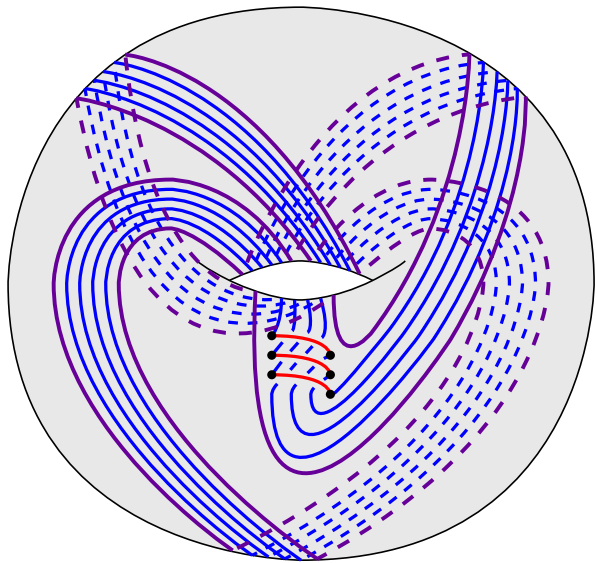}
    \caption{A $(1,3)$-surface $\Sigma_1'$ for $K_1$, with purple curves depicting $T_n \cap \Sigma_1'$}
        \label{bridge1B}
\end{figure}

However, the surface is not minimal; by Theorem \ref{spectra},
\[ b_1(K_1) = \min\{ |21 - 3 \cdot 2 \cdot 4|,4\} = 3.\]
Thus, there is a $(1,3)$-surface $\Sigma_1'$ for $K_1$, shown in Figure \ref{bridge1B}, where $X \cap \Sigma_1'$ is a boundary parallel annulus.  It is not difficult to observe that $\Sigma_1$ is a perturbation of $\Sigma_1'$. \\

Finally, $K_1$ is isotopic into a genus $g$ Heegaard surface for $S^3$ whenever $g \geq 2$, and we have
\[ \mathbf{b}(K_1) = (8,3,0);\]
hence, the bridge spectrum of $K_1$ has two gaps.

\section{Questions}\label{question}

We conclude with several open questions that may be of interest.

\begin{question}
What other spectra can be realized by knots in $S^3$?  Specifically, for any decreasing sequence $\mathbf{v}$ of of positive integers, is there a knot $K$ such that $\mathbf{b}(K) = \mathbf{v}$?
\end{question}

\begin{question}
What other families of knots have interesting bridge spectra?  For instance, what is the bridge spectrum of a twisted torus knot?
\end{question}

While it is relatively straightforward to exhibit candidate bridge surfaces for twisted torus knots, showing these positions to be minimal is a more complicated task.  In general, Lemma \ref{L10} does not apply to hyperbolic twisted torus knots, and so a new strategy would likely be required. \\

The next question is posited in \cite{doll} by Doll:

\begin{question}
For knots $K_1$ and $K_2$ in $S^3$, is there a relationship between $\mathbf{b}(K_1)$, $\mathbf{b}(K_2)$, and $\mathbf{b}(K_1 \# K_2)$?
\end{question}

A simple construction shows that
\[ b_g(K_1 \# K_2) \leq \min \{b_{g_1}(K_1) + b_{g_2}(K_2) - 1: g_1 + g_2 = g\}.\]
However, for any $n$, there are knots $K_1$ and $K_2$ such that $t(K_1 \# K_2) < t(K_1) + t(K_2) - n$ (see \cite{koba} and \cite{morim2}); hence, these knots have the property that for some $g$ the inequality above is strict.  It may be possible that the inequalities above become equalities when we restrict to the class of meridionally small knots, and Lemma \ref{L10} may be of use here. \\

We can also examine the overall bridge structure of iterated torus knots.

\begin{question}
Is there an iterated torus knot $K_n$ with an irreducible $(g,b)$-bridge surface $\Sigma$ such that $b > b_g(K_n)$?
\end{question}

In \cite{schartom}, it is shown that every bridge surface for a 2-bridge knot is the result of stabilization, perturbation, and meridional stabilization performed on a $(0,2)$-surface.  Is it possible that all bridge surfaces for iterated torus knots are derived in this way from the bridge surfaces exhibited here, or are there any unexpected surfaces?

\end{document}